\documentclass[12pt]{amsart} 
\usepackage{ifthen} 
\usepackage{amsmath,amssymb,latexsym} 
\usepackage{times}
\usepackage[T1]{fontenc}

\font\goth=eusm10
\newcommand\F{\mathcal F}
\newcommand\E{\mathcal E}
\newcommand\N{\mathcal N}
\newcommand\T{\mathcal T}
\newcommand\La{\mathcal L}
\newcommand\G{\mathcal G}
\newcommand\V{\mathcal V}

\newcommand\Ss{\mathcal S}
\newcommand\Pmc{\mathcal P}
\newcommand\Ii{\hbox{\goth I}}
\newcommand\Oc{\hbox{\goth O}}
\newcommand\mop{\underline{m}_p}

\newcommand\ZZ{\mathbb{Z}}
\newcommand\CC{\mathbb{C}}

\newcommand\RR{\mathbb{R}}
\newcommand\QQ{\mathbb{Q}}

\newcommand\Pp{\mathbb P}

\newcommand\Pd{\mathbb{P}^2}
\newcommand\Pt{\mathbb{P}^3}

\newcommand\Ct{\tilde{C}}

\numberwithin{equation}{section}

\newtheorem{theorem}[equation]{Theorem}
\newtheorem{coroll}[equation]{Corollary}
\newtheorem{corollary}[equation]{Corollary}

\newtheorem{proposition}[equation]{Proposition}

\newtheorem{definition}[equation]{Definition}

\newtheorem*{theorema}{Theorem}

\newtheorem{remark}[equation]{Remark}

\newtheorem{assumption}{Assumption}

\begin{document} 

\title{Equivalence of families of singular schemes on threefolds and on ruled fourfolds}
\author{Flaminio Flamini}

\email{flamini@matrm3.mat.uniroma3.it}
\curraddr{Dipartimento di Matematica, Università degli Studi de L'Aquila, Via 
Vetoio - Coppito 1, 67100 L'Aquila - Italy}

\thanks{The author is a member of Cofin GVA, EAGER and GNSAGA-INdAM}

\begin{abstract}{The main purpose of this paper is twofold. We first 
want to analyze in details the meaningful geometric aspect of the method 
introduced in the previous paper \cite{F2}, concerning regularity of 
families of irreducible, nodal "curves" on a smooth, projective threefold $X$.  
This analysis highlights several fascinating 
connections with families of other singular geometric 
"objects" related to $X$ and to other varieties. 

Then, we generalize this method to study similar 
problems for families of singular divisors on ruled fourfolds suitably 
related to $X$.}
\end{abstract}

\maketitle

\section*{Introduction}\label{S:intro}

The theory of families of singular curves with fixed 
invariants (e.g. geometric genus, singularity type, 
number of irreducible components, etc.) and which are contained in a 
projective variety $X$ has been extensively 
studied from the beginning of Algebraic Geometry and it 
actually receives a lot of attention, partially due to 
its connections with several fields in Geometry and Physics.

Nodal curves play a central role in the subject of singular curves. 
Families of irreducible and $\delta$-nodal curves on a given projective variety 
$X$ are usually 
called {\em Severi varieties} of irreducible, $\delta$-nodal curves in $X$. 
The terminology "Severi variety" is due to the classical case of families of 
nodal curves on $X= \Pd$, which was first studied by Severi (see \cite{Sev}).

The case in which $X$ is a smooth projective surface has recently given rise 
to a huge amount of literature (see, for example, \cite{CH}, \cite{CC}, \cite{CL}, 
\cite{CS}, \cite{F1}, \cite{GLS}, \cite{Harris}, \cite{Ran}, \cite{S} just to 
mention a few. For a chronological 
overview, the reader is referred for example to Section 2.3 in \cite{F} and 
to its bibliography). This depends not only on the great interest in the subject, but 
also because for a Severi 
variety $V$ on an arbitrary projective variety $X$ there are several problems concerning $V$ 
like non-emptyness, smoothness, irreducibility, dimensional computation as well as 
enumerative and moduli properties of the family of curves 
it parametrizes.

On the contrary, in higher dimension only few results are known. Therefore, in \cite{F2} we focused on 
what is the next relevant case, from the point of view 
of Algebraic Geometry: families of nodal curves on smooth, projective threefolds. 

The aim of this paper is twofold: first, we want to study in details 
the meaningful geometric aspect of the method introduced in \cite{F2}. As a result of this analysis, 
we discover several intriguing and fascinating 
connections with families of other singular geometric 
"objects" related to $X$. Then, we generalize this method to study similar 
problems for families of singular divisors on ruled fourfolds suitably 
related to $X$. 

To be more precise, let $X $ be a smooth projective threefold and let 
$\F$ be a rank-two vector bundle on $X$, which is assumed to be 
globally generated with general global section $s$ having its 
zero-locus $V(s)$ a smooth, irreducible curve $D = D_s$ in 
$X$. The geometric genus of $D$ is given by$$2g(D) -2 = 2p_a(D) -2= deg(\La \otimes \omega_X \otimes 
{\Oc}_D),$$where $\La := c_1(\F) \in Pic(X)$ and $\omega_X$ is the canonical sheaf of $X$.

Take now ${\Pp}(H^0(X, {\F}))$; from our assumptions on $\F$, 
its general point parametrizes a global section whose zero-locus is a 
smooth, irreducible curve.  This projective space  somehow gives a scheme 
dominating a subvariety in which the curves move. 

Given a positive integer $\delta \leq p_a(D)$, 
it makes sense to consider the locally closed
subscheme:
\[
\begin{aligned}
{\V}_{\delta}({\F}) := & \{[s] \in {\Pp}(H^0(X,{\F})) \; | \; C_s := V(s) \subset X \; 
{\rm is \; irreducible} \\
 & {\rm  with \; only} \;  \delta \; {\rm nodes \; as \; singularities} \}; 
\end{aligned}
\](cf. \eqref{eq:severi var}). These are usually called {\em Severi varieties} of global sections of $\F$ 
whose zero-loci are irreducible, $\delta$-nodal curves in $X$, of arithmetic genus $p_a(D)$ 
and geometric genus $g= p_a(D) - \delta$ (cf. \cite{BC}, for $X=\Pt$, and 
\cite{F2} in general). This is because such schemes are 
the natural generalization of the (classical) Severi varieties 
on smooth, projective surfaces recalled before.

When ${\V}_{\delta} ({\F})$ is not empty then
its expected codimension in ${\Pp}(H^0(X,{\F}))$
is $\delta$ (see Proposition \ref{prop:0}). Thus, one says that 
a point $[s] \in {\V}_{\delta} ({\F})$ is a {\em regular 
point} if it is smooth and such that $dim_{[s]}({\V}_{\delta} ({\F}))$ 
equals the expected one (cf. Definition \ref{def:0}). 

It is clear from the definition of regularity that it is fundamental 
to determine the tangent space to a Severi variety at a given point. 
In \cite{F2} we introduced the following cohomological description of 
the tangent space $T_{[s]}({\V}_{\delta} ({\F}))$.

\begin{theorema} {\bf 1} (cf. Theorem \ref{prop:3.fundamental})
Let $X $ be a smooth projective threefold. Let $\F$ be a globally generated rank-two vector
bundle on $X$ and let $\delta$ be a positive integer. Fix $[s] \in {\V}_{\delta}({\F})$ and
let $C=V(s) \subset X$. Denote by $\Sigma$ the set of nodes of
$C$. 

\noindent
Let$${\mathcal P} := {\Pp}_{X}({\F}) \stackrel{\pi}{\longrightarrow} X$$be
the projective space bundle together with its natural projection $\pi$ on $X$ and denote
by ${\Oc}_{\mathcal P}(1)$ its tautological line bundle. 

Then, there exists a zero-dimensional subscheme $\Sigma^1 \subset 
\Pmc$ of length $\delta$, which is a set of $\delta$ rational double points for the
divisor $G_s \in |{\Oc}_{\mathcal P}(1)|$ corresponding to the given section
$s \in H^0(X, \F)$.

\noindent
In particular, each element $[s] \in {\V}_{\delta}({\F}) $ corresponds 
to a divisor $G_s \in |{\Oc}_{\mathcal P}(1)|$, which contains the 
$\delta$ fibres $L_{p_i} = \pi^{-1}(p_i) \subset \Pmc$, for $p_i \in \Sigma$, 
$1 \leq i \leq \delta$, 
and which has $\delta$ rational double points each of which are  
on exactly one of the $\delta$ fibres $L_{p_i}$. 
\end{theorema}

By using the above result, one can translate 
the regularity property of $[s] \in {\V}_{\delta}({\F})$ 
into the surjectivity of some maps among spaces of sections of suitable sheaves
on $X$ (cf. Proposition \ref{prop:corr1} and Remark \ref{rem:localdescr}). 
This allows us to find several equivalent and sufficient 
conditions for the regularity of Severi varieties 
${\V}_{\delta}({\F})$ on $X$ (cf. \cite{F2} and also 
Corollaries \ref{cor:propcorr1}, \ref{rem:reg}, \ref{cor:explanetion}, 
\ref{cor:tangentspace} and Theorem \ref{thm:9bis} in this paper).

On the other hand, Theorem 1 also introduces a meaningful and fascinating 
connection between elements in 
$ {\V}_{\delta}({\F}) $ and other singular schemes in $X$ and in $\Pmc$. 
The aim of this paper is to investigate in details 
the deep geometric meaning of $T_{[s]}({\V}_{\delta}({\F}))$ and its 
several connections with these 
families of singular geometric objects related to $X$ and to $\Pmc$. 

More precisely, we show that there exists a natural connection among: 
\begin{itemize}
\item a nodal section $[s] \in {\V}_{\delta}({\F})$ on $X$,
\item the corresponding singular divisor $G_s$ in $|\Oc_{\Pmc}(1) |$ on $\Pmc$ having $\delta$ 
rational double points over the nodes of $C = V(s)$, 
\item for any $s + \epsilon \; s' \in T_{[s]}({\V}_{\delta}({\F}))$, 
the surface $V(s \wedge s') \subset X$ which is singular along $\Sigma$ and which belongs 
to the linear system $| \Ii_{C/X} \otimes \La|$ on $X$, 
\item given $G_s$ and $G_{s'}$ in $|\Oc_{\Pmc}(1) |$ corresponding to $s$ and $s'$ respectively, 
the "complete intersection" surface in $\Pmc$, 
${\Ss}_{s,s'} := G_s \cap G_{s'}$, which is singular along $\Sigma^1$ and which dominates 
$V(s \wedge s')$, 
\end{itemize} (cf. Propositions \ref{prop:corr1}, \ref{prop:corr2} and 
Remarks \ref{rem:propcorr2}, \ref{rem:2propcorr2}). 

In particular, we prove:

\begin{theorema} {\bf 2} (cf. Theorem \ref{thm:fundamental} and Corollary \ref{cor:explanetion}) 
The following conditions are equivalent:
\begin{itemize}
\item[(i)] $s + \epsilon s' \in T_{[s]}({\V}_{\delta}({\F}))$, where $\epsilon^2 = 0$;
\item[(ii)] $V (s \wedge s') \subset X$ is a surface which 
contains $C$ and which is singular along $\Sigma$;
\item[(iii)] the divisor $G_{s'}$ passes through $\Sigma^1$
\item[(iv)] the surface $\Ss_{s,s'}:= G_s \cap G_{s'} \subset \Pmc$ is singular 
along $\Sigma^1$; 
\end{itemize} 
\end{theorema}By using Theorem 2 we can also give several further 
equivalent conditions for the 
regularity of $[s] \in {\V}_{\delta}({\F})$ (cf. Remark \ref{rem:explanetion}).

Furthermore, thanks to the correspondence introduced in Theorems 1 and 2, we also consider 
a generalization of Severi varieties of nodal curves. Indeed, we denote by 
$${\mathcal R}_{\delta} (\Oc_{\Pmc}(1)) 
:=  \{G_s \in |\Oc_{\Pmc}(1) | \; {\rm s.t.} \; [s] \in {\V}_{\delta}({\F}) \}$$the schemes 
parametrizing families of expected codimension $\delta$ 
in $|\Oc_{\Pmc}(1) |$, whose elements correspond to divisors 
which are irreducible and with only $\delta$ rational double points as singularities. 
For brevity sake, these are 
called $\Pmc_{\delta}$-{\em Severi varieties} (cf. Definition \ref{def:rdelta} and 
\eqref{eq:expdim}). 

One can obviously give a similar definition of {\em regularity} 
for $\Pmc_{\delta}$-{\em Severi varieties} 
(cf. Definition \ref{def:expdim}). We prove:

\begin{theorema} {\bf 3} (cf. Theorem \ref{thm:tangentspace} and Corollary \ref{cor:tangentspace}) 
Let $[G_s] \in {\mathcal R}_{\delta} (\Oc_{\Pmc}(1))$ on $\Pmc$ and let $\Sigma^1 $ be 
the zero-dimensional scheme of the $\delta$-rational double points of $G_s \subset \Pmc$. Then
$$T_{[G_s]} ({\mathcal R}_{\delta} (\Oc_{\Pmc}(1))) \cong \frac{H^0(\Ii_{\Sigma^1/\Pmc} 
\otimes \Oc_{\Pmc}(1))}{<G_s>}.$$

In particular,
$$[G_s] \in {\mathcal R}_{\delta} (\Oc_{\Pmc}(1)) \; {\rm is \; a \; regular \; point}   \Leftrightarrow 
 [s] \in {\V}_{\delta}({\F}) \; {\rm is \; a \; regular \; point}.$$
\end{theorema}

Finally, we first improve some regularity results of \cite{F2} 
for Severi varieties ${\V}_{\delta}({\F})$ of irreducible, 
$\delta$-nodal sections on $X$; then, we use Theorem 3 
to deduce regularity results also for $\Pmc_{\delta}$-Severi varieties 
${\mathcal R}_{\delta} (\Oc_{\Pmc}(1)) $ on 
$\Pmc$. Precisely, we have:

\begin{theorema} {\bf 4} (cf. Theorems \ref{thm:9bis} and \ref{thm:regrdelta}) 
Let $X$ be a smooth projective threefold, $\E$ be a globally generated rank-two vector bundle on 
$X$, $M$ be a very ample line bundle on $X$ and $k \geq 0$ and $\delta >0$ be integers. 
Let $\Pmc : = \Pp_X (\E \otimes M^{\otimes k})$ and $\Oc_{\Pmc}(1)$ be its tautological 
line bundle. If$$(*)\;\;\;\;\;\;\delta \leq k+1,$$then both 
${\V}_{\delta}({\E} \otimes M^{\otimes k})$ on $X$ and  ${\mathcal R}_{\delta} (\Oc_{\Pmc}(1))$ 
on $\Pmc$ are regular at each point.
\end{theorema} 

\noindent
The upper-bounds in $(*)$ 
are also shown to be almost-sharp (cf. Remark \ref{rem:10}).

What we want to stress is the following fact: the regularity condition 
for the schemes ${\mathcal R}_{\delta} (\Oc_{\Pmc}(1)) $ on $\Pmc$ 
is equivalent to the separation of suitable 
zero-dimensional schemes by the linear 
system $|\Oc_{\Pmc}(1)|$ on the fourfold $\Pmc$ (cf. Corollary \ref{cor:tangentspace}). 
In general, it is well-known how difficult is to enstablish separation 
of points in projective varieties of dimension greater than or equal to three 
(cf. e.g. \cite{AS}, \cite{EL} and \cite{K}). In some cases, some separation 
results can be found by using technical tools 
like {\em multiplier ideals} as well as the Nadel and the Kawamata-Viehweg 
vanishing theorems (see, e.g. \cite{E}, 
for an overview). In our situation, thanks to the correspondence between 
${\V}_{\delta}({\F})$ on $X$ and 
${\mathcal R}_{\delta} (\Oc_{\Pmc}(1))$ on $\Pmc$, we deduce regularity 
conditions for ${\mathcal R}_{\delta} (\Oc_{\Pmc}(1))$ from those already obtained 
for ${\V}_{\delta}({\F})$.

The paper consists of seven sections. Section \ref{S:1} contains some terminology and notation. In 
Section \ref{S:a} we briefly recall some fundamental definitions in \cite{F2}, which 
are frequently used in the whole paper. Section \ref{S:b} briefly recall one of the main result 
in \cite{F2} concerning the correspondence beteween elements in  ${\V}_{\delta}({\F})$ on $X$ and 
those in ${\mathcal R}_{\delta} (\Oc_{\Pmc}(1))$ on $\Pmc$ (cf. Theorem \ref{prop:3.fundamental}). 

In Section \ref{S:c} we describe how to associate elements of $T_{[s]}{\V}_{\delta}({\F})$ 
to singular divisors in $X$. Section \ref{S:d} contains one of the main result of the paper 
(cf. Theorem \ref{thm:fundamental}) which proves the equivalence of several singular geometric 
"objects" related to $X$ and to $\Pmc$. In Section \ref{S:ef} we focus on 
$\Pmc$-Severi varieties on $\Pmc$; we give a description 
of tangent spaces at points of such schemes as well as we find conditions 
for their regularity (cf. Theorem \ref{thm:tangentspace} and Corollary \ref{cor:tangentspace}). 
Section \ref{S:ghil} is devoted to the determination of 
almost-sharp upper-bounds on $\delta$ implying the regularity 
of ${\V}_{\delta} ({\F})$ on $X$ as well as of ${\mathcal R}_{\delta} (\Oc_{\Pmc}(1))$ on 
$\Pmc$.

\vskip 10pt 

\noindent
{\it Acknowledgments:} Part of this paper is related to the previous paper 
\cite{F2}, which was prepared during my permanence 
at the Department of Mathematics of the University of Illinois at Chicago 
(February - May 2001). Therefore, my deepest gratitude goes to L. Ein, not 
only for the organization of my visit, but mainly for all I have learnt from him as well as 
for having suggested me to approach this pioneering area. 
My very special thanks go also 
to GNSAGA-INdAM and to V. Barucci, A. F. Lopez and E. Sernesi for their confidence 
and their support during my period in U.S.A. 

I am greatful to L. Caporaso and to the organizers of the Workshop {\em "Global geometry of 
algebraic varieties"} - Madrid, December 16-19, 2002. In fact, their invitations to give talks 
at the Dept. of Mathematics of "Roma Tre" and at the Workshop, respectively, have given "life" 
to these new results, because of the preparation of the talks on the subjects in \cite{F2}. 
I am indebted to L. Chiantini and to C. Ciliberto for having always stimulated me to 
investigate more in this research area.

\section{Notation and Preliminaries}\label{S:1}
We work in the category of algebraic 
$\CC$-schemes. $Y$ is a \emph{$m$-fold} if it is a reduced, irreducible and non-singular scheme 
of finite type and of dimension $m$. 
If $m=1$, then $Y$ is a (smooth) {\em curve}; $m=2$ and $3$ are the cases of a 
(non-singular) {\em surface} and {\em threefold}, respectively. 
If $Z$ is a closed subscheme of a scheme $Y$, $\Ii_{Z/Y}$ 
denotes the \emph{ideal sheaf} of $Z$ in $Y$, ${\N}_{Z/Y}$ 
the {\em normal sheaf} of $Z$ in $Y$ whereas 
${\N}_{Z/Y}^{\vee} \cong {\Ii_{Z/Y}}/{\Ii_{Z/Y}^2}$ is the {\em 
conormal sheaf} of $Z$ in $Y$. As usual, $h^i(Y, \; -):=\text{dim} \; H^i(Y, \; -)$.

Given $Y$ a projective scheme, $\omega_Y$ denotes its dualizing sheaf. 
When $Y$ is a smooth variety, then $\omega_Y$ coincides with its canonical bundle and 
$K_Y$ denotes a canonical divisor s.t. $\omega_Y \cong \Oc_Y(K_Y)$. Furthermore, 
${\T}_Y$ denotes its tangent bundle whereas $\Omega^1_Y$ denotes its cotangent bundle. 

If $D$ is a reduced curve, $p_a(D)=h^1(\Oc_D)$ 
denotes its {\em arithmetic genus}, 
whereas $g(D)= p_g(D)$ denotes its \emph{geometric genus}, the
arithmetic genus of its normalization.

Let $Y$ be a projective $m$-fold and $\E$ be a rank-$r$ vector bundle on 
$Y$; $c_i(\E)$ denotes the \emph{$i^{th}$-Chern class} of $\E$, 
$1 \leq i \leq r$. As in \cite{Ha} - Sect. II.7 - 
$\Pp_{Y} ({\E})$ denotes the {\it projective space bundle} on $Y$, 
defined as $Proj(Sym({\E}))$. 

There is a surjection 
$\pi^*(\E) \to \Oc_{\Pp_{Y} ({\E})}(1)$, where $\Oc_{\Pp_{Y} ({\E})}(1)$ is the 
{\it tautological line bundle} on $\Pp_{Y} ({\E})$ and where $\pi : \Pp_{Y} ({\E}) \to Y$ is 
the natural projection morphism. Recall that $\E$ is said to be an {\em ample} (resp. 
{\em nef}) vector bundle on $Y$ if $\Oc_{\Pp_{Y} ({\E})}(1)$ is an ample (resp. nef) 
line bundle on $\Pp_{Y} ({\E})$ (see, e.g. \cite{H1}).

For non reminded terminology, the reader is referred to \cite{BPV}, \cite{Fr} and 
\cite{Ha}.

\section{Families of nodal "curves" on smooth, projective threefolds}\label{S:a}

In this section we briefly recall some definitions and 
results from \cite{F2} which will be 
frequently used in the sequel.

Let $X $ be a smooth projective threefold and let $\F$ be a rank-two vector 
bundle on $X$. If $\F$ is  
globally generated on $X$, it is not restrictive if from now on 
we assume that the zero-locus $V(s)$ of its
general global section $s$ is a smooth, irreducible curve $D = D_s$ in 
$X$ (for details, see \cite{F2}; 
for general 
motivations and backgrounds, the reader is referred to e.g. \cite{OSS} and 
to \cite{Sz}, Chapter IV).

From now on, denote by $\La \in Pic(X)$ the 
line bundle on $X$ given by $c_1(\F)$. Thus, by the {\em Koszul sequence} 
of $({\F}, s)$:
\begin{equation}\label{eq:koszul}
0 \to {\Oc}_X \to {\F} \to {\Ii}_{V(s)} \otimes \La \to 0,
\end{equation}we compute 
the geometric genus of $D$ in terms of the invariants of $\F$ and of $X$. Precisely 
\begin{equation}\label{eq:numeriX}
2g(D) -2 = 2p_a(D) -2= deg(\La \otimes \omega_X \otimes 
{\Oc}_D). 
\end{equation}

This integer is easily 
computable when, for example, $X$ is a general complete intersection threefold. 
In particular, when $Pic(X) \cong \ZZ$ (e.g $X= \Pt$ or $X$ either a prime Fano or a 
complete intersection 
Calabi-Yau threefold) one can use this isomorphism to identify line bundles on $X$ 
with integers. 
Therefore, if $A$ denotes the ample generator class of $Pic(X)$ over $\ZZ$ and 
if $\F$ is a rank-two vector bundle on $X$ such that $c_1({\F}) = n A$, we can also 
write $c_1({\F}) = n$ with no ambiguity. 

Thus, if e.g. $X = \Pt$ and if 
we put $c_i = c_i({\F}) \in \ZZ$, we have
\begin{equation}\label{eq:numeriP3}
deg(D) = c_2 \; {\rm and} \; g(D)= p_a(D) = \frac{1}{2} (c_2 (c_1 -4))+ 1, 
\end{equation}i.e. $D$ is subcanonical of level $(c_1 -4)$.

Take now ${\Pp}(H^0(X, {\F}))$; from our assumptions on $\F$, 
the general point of this projective space 
parametrizes a global section whose zero-locus is a 
smooth, irreducible curve in $X$. Given a positive integer $\delta \leq p_a(D)$, 
it makes sense to consider the subset
\begin{equation}\label{eq:severi var}
\begin{aligned}
{\V}_{\delta}({\F}) := & \{[s] \in {\Pp}(H^0(X,{\F})) \; | \; C_s := V(s) \subset X \; 
{\rm is \; irreducible} \\
 & {\rm  with \; only} \;  \delta \; {\rm nodes \; as \; singularities} \}; 
\end{aligned}
\end{equation}therefore, any element of ${\V}_{\delta}({\F})$ determines a curve in $X$ 
whose arithmetic genus $p_a(C_s)$ is given by 
(\ref{eq:numeriX}) and whose geometric genus is $g= p_a(C_s) - \delta$. 
We recall that ${\V}_{\delta}({\F})$ is a locally closed
subscheme of the projective space ${\Pp}(H^0(X, {\F}))$; it 
is usually called the {\em Severi variety} of global sections of $\F$ 
whose zero-loci are irreducible, 
$\delta$-nodal curves in $X$ (cf. \cite{BC}, for $X=\Pt$, and 
\cite{F2} in general). This is because such schemes are 
the natural generalization of the (classical) Severi varieties of irreducible and $\delta$-nodal curves 
in linear systems on smooth, projective surfaces (see \cite{CC}, \cite{CH},
\cite{CS}, \cite{F1}, \cite{GLS}, \cite{Harris}, \cite{Ran}, \cite{S} and \cite{Sev}, just to 
mention a few). 

For brevity sake, we shall usually refer to ${\V}_{\delta}({\F})$ as the 
{\em Severi variety of irreducible, $\delta$-nodal sections} of $\F$ on $X$.

First possible questions on such Severi varieties are about 
their dimensions as well as their smoothness properties. 

A preliminary estimate is given by the following standard result:

\begin{proposition}\label{prop:0}
Let $X $ be a smooth projective threefold, $\F$ be a 
globally generated rank-two vector bundle on $X$ and $\delta$ be a positive integer. 
Then
\[expdim ({\V}_{\delta}({\F})) = \begin{cases}
				 h^0(X, {\F}) - 1 - \delta, \; \; {\rm if} \; \delta \leq 
                                 h^0(X, {\F}) - 1= dim({\Pp}(H^0({\F}))),\\ 
 				 -1, \;\;\;\;\;\;\;\;\;\;\;\;\;\;\;\;\;\;\;\;\;\;\;\;{\rm if} 
                                 \; \delta \geq 
                                 h^0(X, {\F}).\\ 
				\end{cases} \]

\end{proposition}
\begin{proof} See Proposition 2.10 in \cite{F2}.
\end{proof}

\begin{assumption}\label{ass:main}
\normalfont{From now on, 
given $X$ and $\F$ as in Proposition \ref{prop:0},
we shall always assume ${\V}_{\delta}({\F}) \neq \emptyset$. We write $[s] \in {\V}_{\delta}({\F})$ 
to intend that the global section $s \in H^0(X, {\F})$ determines the corresponding point $[s]$ 
of the scheme ${\V}_{\delta}({\F})$. We simply denote by $C$ - instead of $C_s$ - the zero-locus of the given 
section $s$, when it is clear from the context that $s$ is fixed. We finally consider 
$\delta \leq {\rm min} \{ h^0(X, {\F}) - 1, \; p_a(C) \}$ - the latter 
is because we want $C=V(s)$ to be irreducible, for any $[s] \in {\V}_{\delta}({\F})$.
}
\end{assumption}

By Proposition \ref{prop:0}, it is natural to state the following:

\begin{definition}\label{def:0}
Let $[s] \in {\V}_{\delta}({\F})$, with 
$\delta \leq {\rm min} \{ h^0(X, {\F}) - 1, \; p_a(C) \}$.
Then $[s]$ is said to be a {\em regular point} of ${\V}_{\delta}({\F})$ if:
\begin{itemize}
\item[(i)] $[s] \in {\V}_{\delta}({\F})$ is a smooth point, and
\item[(ii)] $dim_{[s]}({\V}_{\delta}({\F})) = expdim ({\V}_{\delta}({\F}))=
dim({\Pp}(H^0(X,{\F}))) - \delta$.
\end{itemize}${\V}_{\delta}({\F})$ is said to be {\em regular} if it is 
regular at each point. 
\end{definition}

One of the main result in \cite{F2} has been to
present a cohomological description
of the tangent space $T_{[s]}({\V}_{\delta}({\F}))$ which translates
the regularity property of a given point $[s] \in {\V}_{\delta}({\F})$ into the
surjectivity of some maps among spaces of sections of suitable sheaves
on the threefold $X$. This description allowed us to find also several 
sufficient conditions for the regularity of Severi varieties 
${\V}_{\delta}({\F})$ on $X$.

One of the aim of this paper is to study in more details 
the deep geometric meaning of the cohomological description of 
the tangent space $T_{[s]}({\V}_{\delta}({\F}))$ and its 
several connections 
with families of other singular geometric objects related 
to $X$ and to $\F$.

To do this, we have first to recall some results contained in \cite{F2}, 
since these are the starting point of our analysis.

\section{Association of elements of ${\V}_{\delta}({\F})$ on $X$ to 
singular divisors in  ${\Pp}_{X}({\F}) $}\label{S:b}

In this section we want to briefly recall the 
correspondence given in \cite{F2} between elements of ${\V}_{\delta}({\F})$ on $X$  
and suitable singular divisors in the tautological linear system $|{\Oc}_{\mathcal P}(1)|$ 
on the projective space bundle 
$\Pmc :={\Pp}_{X}({\F})$, which is a fourfolds 
ruled over $X$. Instead of referring the reader to \cite{F2}, we prefer to briefly 
recall here the proofs of some results contained in there, not only because  
we give here more precise statements and proofs, 
but mainly because the strategy of the proofs as well as their technical details will be 
fundamental for the analysis in the whole paper.

From now on, with conditions as in Assumption \ref{ass:main}, 
let $[s] \in {\V}_{\delta}({\F})$. 
Then:

\begin{theorem}\label{prop:3.fundamental} (cf. Theorem 3.4 (i) in \cite{F2})
Let $X $ be a smooth projective threefold. Let $\F$ be a globally generated rank-two vector
bundle on $X$ and let $\delta$ be a positive integer. Fix $[s] \in {\V}_{\delta}({\F})$ and
let $C=V(s) \subset X$. Denote by $\Sigma$ the set of nodes of
$C$. 

\noindent
Let$${\mathcal P} := {\Pp}_{X}({\F}) \stackrel{\pi}{\longrightarrow} X$$be
the projective space bundle together with its natural projection $\pi$ on $X$ and denote
by ${\Oc}_{\mathcal P}(1)$ its tautological line bundle. 

Then, there exists a zero-dimensional subscheme $\Sigma^1 \subset 
\Pmc$ of length $\delta$, which is a set of $\delta$ rational double points for the
divisor $G_s \in |{\Oc}_{\mathcal P}(1)|$ corresponding to the given section
$s \in H^0(X, \F)$.

\noindent
In particular, each element $[s] \in {\V}_{\delta}({\F}) $ corresponds 
to a divisor $G_s \in |{\Oc}_{\mathcal P}(1)|$, which contains the 
$\delta$ fibres $L_{p_i} = \pi^{-1}(p_i) \subset \Pmc$, for $p_i \in \Sigma$, 
$1 \leq i \leq \delta$, 
and which has $\delta$ rational double points each of which are  
on exactly one of the $\delta$ fibres $L_{p_i}$. 
\end{theorem}
\begin{proof} 
Consider the smooth, projective, ruled fourfold
$${\mathcal P} := {\Pp}_{X}({\F}) \stackrel{\pi}{\longrightarrow} X,$$together
with its tautological line bundle ${\Oc}_{\mathcal P}(1)$ such that
$\pi_*({\Oc}_{\mathcal P}(1)) \cong {\F}$. Since $[s] \in {\V}_{\delta}({\F}) $, in particular 
$s \in H^0 (X, \F)$; then, one also has
$$0 \to {\Oc}_{\Pmc} \stackrel{\cdot s}{\to} {\Oc}_{\Pmc}(1).$$Therefore, the nodal curve
$C \subset X$ corresponds to a divisor - say $G_s $ - on the
fourfold $\Pmc$ which belongs to the tautological linear system $|{\Oc}_{\Pmc}(1) |$. 

It is clear that $G_s$ contains all the $\pi$-fibres over 
the zero-locus $C=V(s)$; precisely, it contains the surface 
$\mathbb{F}:= {\Pp}^1_C = Proj(Sym(\F|_C))$ which is ruled over $C$. Therefore, 
in particular $G_s$ contains the locus $\Lambda := \bigcup_{i=1}^{\delta} L_{p_i}$, 
where $\Sigma = \{p_1, \ldots, p_{\delta} \}$. 

The geometry of $G_s$ is strictly related 
to the one of $C$. Indeed, if $p \in \Sigma = Sing(C)$, take $U_p \subset X$
an affine open set containing $p$, where the vector bundle ${\F}$ trivializes. 
Since $p$ is a {\em planar singularity}, one can choose local coordinates 
$\underline{x} = (x_1, x_2, x_3)$ on $U_p \cong {\mathbb A}^3$ such
that $\underline{x}(p) = (0,0,0)$ and such that the global section $s$ is
$$s|_{U_p} = (x_1x_2, \; x_3).$$For 
what concerns the divisor 
$G_s \in |{\Oc}_{\mathcal P}(1)|$ corresponding to $s \in H^0 (X, \F)$, since
$U_p$ trivializes $\F$, then ${\Pmc} \; |_{U_p} \cong U_p \times {\Pp}^1.$ Taking homogeneous
coordinates $[u,v] \in {\Pp}^1$,we have ${\Oc}_{\Pmc}(\pi^{-1}(U_p)) \cong 
{\CC}[x_1, x_2, x_3, u, v]$. Thus,
$$ {\Oc}_{G_s}(\pi^{-1}(U_p)) \cong {\CC}[x_1, x_2, x_3, u, v]/(u x_1 x_2 + v x_3).$$This implies 
that the local equation of $G_s$ in $\pi^{-1}(U_p)$ is given by
\begin{equation}\label{eq:loccomp0}
u x_1 x_2 + v x_3 = 0. 
\end{equation}

In the open chart where $v \neq 0$, $G_s$ is smooth, whereas where 
$u \neq 0$, we see that \eqref{eq:loccomp0} 
is the equation of a quadric cone in $\mathbb{A}^4$ having vertex at the origin. This means that 
$G_s$ has a rational double point along the $\pi$-fibre $L_p:= \pi^{-1} (p) \subset \Pmc$.

Globally speaking, one can state that there exist 
$\delta$ distinguished points on $\Pmc$. Such 
distinguished points determine a $0$-dimensional subscheme 
$$\Sigma^1 \subset {\Pmc}$$along which the divisor $G_s \subset \Pmc$ is singular. Thus 
$\Sigma^1$ is a set of 
$\delta$ rational double points for $G_s$, each line of 
${\Lambda} = \pi^{-1}(\Sigma)= \bigcup_{i=1}^{\delta} L_{p_i} $ 
containing only one of such $\delta$ points; more precisely, each point 
$p^1_i\in \Sigma^1$ is a rational double point for $G_s$ and it belongs 
to the fibre $L_{p_i} \subset \Pmc$, where $p_i \in \Sigma$ is a node of 
$C = V(s)$. Therefore, the isomorphism$$ \Sigma^1 \cong \Sigma$$is directly given 
by the natural projection $\pi : \Pmc \to X$.
\end{proof}

The above result introduces a 
correspondence between elements of ${\V}_{\delta}({\F})$ and suitable 
divisors in $|{\Oc}_{\mathcal P}(1)|$; this correspondence has been used 
in \cite{F2} to give a cohomological description of the tangent space 
$T_{[s]}({\V}_{\delta}({\F}))$, which has been a fundamental point 
in order to determine several sufficient 
conditions for the regularity of Severi varieties ${\V}_{\delta}({\F})$ on 
$X$ (cf. Theorems 4.5, 5.9, 5.25, 5.28 and 5.36 in \cite{F2}).  

The main ideas are as follows: $C$ is local complete intersection in $X$, whose 
normal sheaf is the rank-two vector-bundle ${\N}_{C/X} \cong {\F}|_C$. 
Let $T^1_C$ be the {\em first cotangent sheaf} of $C$, i.e.
$T^1_C \cong {\mathcal Ext}^1(\Omega^1_C, {\Oc}_C)$,
where $\Omega_C^1$ is the sheaf of K$\ddot{a}$hler differentials of the nodal
curve $C$ (for details, see \cite{LS}). Since $C$ is nodal, 
$T^1_C$ is a sky-scraper sheaf supported on $\Sigma$, such that 
$T^1_C \cong \bigoplus_{i=1}^{\delta} {\CC}_{(i)}$. Furthermore, one 
has the exact sequence:
\begin{equation}\label{eq:T1}
0 \to {\N}'_C \to {\N}_{C/X} \stackrel{\gamma}{\to} T^1_C \to 0, 
\end{equation}where ${\N}'_C$ is defined as the kernel of the natural surjection 
$\gamma$ (see, for example, \cite{S}).

\begin{proposition}\label{prop:corr1} (cf. Theorem 3.4 (ii) in \cite{F2}) 
With assumptions and notation as in Theorem \ref{prop:3.fundamental}, denote 
by $\Ii_{\Sigma^1/{\Pmc}}$ the ideal sheaf of $\Sigma^1$ in $\Pmc$. 
Then the subsheaf of $\F$, defined by
\begin{equation}\label{eq:svolta}
{\F}^{\Sigma} := \pi_* ({\Ii}_{\Sigma^1/{\Pmc}} \otimes {\Oc}_{\mathcal P}(1)),
\end{equation}is such
that its global sections (modulo the
one dimensional subspace $<s>$) parametrize first-order deformations of
$s \in H^0(X, \F)$ which are equisingular. 

Precisely, we have
\begin{equation}\label{eq:reg}
\frac{H^0(X, {\F}^{\Sigma})}{< s >} \cong T_{[s]} ({\V}_{\delta}({\F})) \subset T_{[s]}
({\Pp}(H^0({\F}))) \cong \frac{H^0(X, {\F})}{< s >}.
\end{equation}
\end{proposition}
\begin{proof} By the correspondence given 
in Theorem \ref{prop:3.fundamental}, we can consider the closed immersion 
$ \Sigma^1 \subset \Pmc$ and so the natural exact sequence
\begin{equation}\label{eq:idealP}
0 \to {\Ii}_{\Sigma^1/{\Pmc}} \otimes {\Oc}_{\mathcal P}(1) \to {\Oc}_{\mathcal P}(1)
\to {\Oc}_{\Sigma^1} \to 0,
\end{equation}which is defined by restricting ${\Oc}_{\mathcal P}(1) $ to $\Sigma^1$. 

Since $\pi_*({\Oc}_{\mathcal P}(1)) \cong {\F}$,  
$\pi_*({\Oc}_{\Sigma^1}) = \pi_*({\Oc}_{\pi^{-1}(\Sigma^1)}) = \pi_*(\pi^*({\Oc}_{\Sigma}))
\cong {\Oc}_{\Sigma}$ 
and since we have ${\F} \to \!\!\! \to {\Oc}_{\Sigma}$, by applying $\pi_*$ to the exact sequence
(\ref{eq:idealP}), we get ${\mathcal R}^1 \pi_*({\Ii}_{\Sigma^1/{\Pmc}}
\otimes {\Oc}_{\mathcal P}(1)) = 0$. Thus, we define
${\F}^{\Sigma}$ as in \eqref{eq:svolta}, so that
\begin{equation}\label{eq:idealX}
0 \to   {\F}^{\Sigma}  \to {\F} \to   {\Oc}_{\Sigma} \to 0,
\end{equation}holds.

Observe that ${\F}^{\Sigma} $ fits in the following exact diagram:
\begin{equation}\label{eq:(*1)} 
\begin{aligned}
\begin{array}{rcccclr}
 &0 & & 0 & & & \\
 & \downarrow & & \downarrow & & & \\
0 \to   &{\Ii}_{C/X} \otimes {\F}  & \stackrel{\cong}{\to} & {\Ii}_{C/X}\otimes {\F} &\to & 0
&  \\ 
 & \downarrow & & \downarrow & & \downarrow & \\
 0 \to  & {\F}^{\Sigma} & \to & {\F} & \to  & \Oc_{\Sigma} & \to 0 \\ 
 & \downarrow & & \downarrow & & \downarrow^{\cong}& \\
0 \to & {\N}'_C & \to &  {\F} |_{C} & \to & T^1_C & \to 0  \\ 
 &\downarrow & & \downarrow & &\downarrow & \\
 & 0 & & 0 & & 0 & .
\end{array}
\end{aligned}
\end{equation}

From the commutativity of diagram \eqref{eq:(*1)}, the vector space 
$$\frac{H^0(X, {\F}^{\Sigma})}{< s >}$$parametrizes the first-order 
deformations of $[s]$ in $\Pp(H^0(X, {\F})) $ which are {\em equisingular}; 
indeed, these are exactly the global sections of ${\F}$ which go to zero at $\Sigma$ in the 
composition 
\begin{equation}\label{eq:compo}
{\F} \to \!\! \to {\F}|_C \to \!\! \to T^1_C \cong {\Oc}_{\Sigma}.
\end{equation}
\end{proof}

Notice that \eqref{eq:reg} gives a completely general characterization of the tangent 
space $T_{[s]} ({\V}_{\delta}({\F}))$ on $X$. In particular, with assumptions and notation as in 
Theorem \ref{prop:3.fundamental} and in Proposition \ref{prop:corr1}, 
we get the following results:

\begin{coroll}\label{cor:propcorr1}
Each global section $s' \in H^0(X, {\F}^{\Sigma})$ corresponds to a 
divisor $G_{s'} \in | \Ii_{\Sigma^1/ \Pmc} 
\otimes \Oc_{\Pmc}(1) |$ on $\Pmc$. In particular, for $\epsilon \in \CC[T]/(T^2)$, we have:

\vskip5pt

$s + \epsilon s' \in T_{[s]} ({\V}_{\delta}({\F}))$ on $X$ $\Leftrightarrow$ $s' \in 
H^0(X, {\F}^{\Sigma})$ $\Leftrightarrow$ $G_{s'} \in |\Oc_{\Pmc}(1) |$ and $\Sigma^1 \subset 
G_{s'}$.
\end{coroll}
\begin{proof}
It directly follows from the definition of $\F^{\Sigma}$ and from the correspondence 
in Theorem \ref{prop:3.fundamental}.
\end{proof}

\begin{coroll}\label{rem:reg} (cf. Corollary 3.9 in \cite{F2}) 
From \eqref{eq:idealX}, it follows that
\begin{equation}\label{eq:regbis}
\begin{array}{rcl}
[s] \in {\V}_{\delta}({\F}) \; {\rm is \; regular} & \Leftrightarrow &
H^0(X , {\F}) \stackrel{\mu_X}{\to \!\!\! \to} H^0(X, {\Oc}_{\Sigma})\\
 & \Leftrightarrow & 
H^0({\Pmc} , {\Oc}_{\mathcal P}(1)) \stackrel{\rho_{\Pmc}}{\to \!\!\! \to} H^0({\Pmc},{\Oc}_{\Sigma^1}).
\end{array}
\end{equation}
\end{coroll}
\begin{proof}
It follows from Proposition \ref{prop:0}, from Theorem \ref{prop:3.fundamental} and from Proposition 
\ref{prop:corr1}.
\end{proof}

\begin{remark}\label{rem:localdescr}
\normalfont{
Note that, on the one hand, the map $\mu_X$ in \eqref{eq:regbis} is not defined by restricting 
the global sections of $\F$ to $\Sigma$ because (\ref{eq:idealX}) - i.e.
the second row of diagram \eqref{eq:(*1)} -
does not coincide with the restriction sequence
$$0 \to   {\Ii}_{\Sigma/X} \otimes {\F}  \to {\F} \to   {\F} |_{\Sigma} \to 0;$$indeed 
${\F} |_{\Sigma}$ has rank two at each node, whereas $\F^{\Sigma}$ has rank one at each node.

On the other hand, by the Leray isomorphism the exact sequence 
(\ref{eq:idealP}) on the fourfold $\Pmc$ is equivalent in cohomology 
to the one in (\ref{eq:idealX}) but it is more naturally defined
by restricting the line bundle ${\Oc}_{\mathcal P}(1)$ to $\Sigma^1$. Therefore, the map 
$\rho_{\Pmc}$ in (\ref{eq:regbis}) is a standard restriction map.

To better understand the geometric meaning of the 
map $\mu_X$, we briefly recall the local description of 
(\ref{eq:idealX}). Therefore, it suffices to consider $\delta =1$.  Assume $\{p \} = Sing(C) = \Sigma$ 
and take, as before, $U_p \subset X$ an affine open set containing $p$, where the vector bundle 
${\F}$ is trivial. Take local 
coordinates $\underline{x} = (x_1, x_2, x_3)$ 
on $U_p \cong {\mathbb A}^3 = \CC^3$ such 
that $\underline{x}(p) = (0,0,0)$ and such that the global section $s$, whose zero-locus 
is $C$, is given by $s|_{U_p} = (x_1x_2, \; x_3).$ Since $C = V(x_1x_2, \; x_3) \subset  U_p$, around the node 
$ \underline{x}(p)= \underline{0} $ the Jacobian map:
$$(**) \;\;\;\;  {\T}_{{\CC}^3}|_C \stackrel{J(s)}{\longrightarrow} {\N}_{C/{\CC}^3} \to T^1_C $$is given 
by  
\[J(s) := \left( \begin{array}{ccc}
			\frac{\partial f_1}{\partial x_1} & \frac{\partial f_1}{\partial x_2} & 
                         \frac{\partial f_1}{\partial x_3}\\
                         \frac{\partial f_2}{\partial x_1} & \frac{\partial f_2}{\partial x_2} & 
                         \frac{\partial f_2}{\partial x_3}
		 	\end{array}
			\right) = 
                    \left( \begin{array}{ccc}
			x_2 & x_1 & 
                         0 \\
                       0  & 0  & 1
		 	\end{array}
			\right). \]Put $Im(J(s)) = < x_2 e'_1, x_1 e'_1, e'_2>$, where $\{e'_1, e'_2 \}$ 
a local basis for ${\N}_{C/{\CC}^3}$. Thus, $e'_2$ goes to zero in $T^1_C$ so, 
by this local description, it follows that the map $\mu_X$ is exactly the composition of 
the evaluation at $p$ of global sections together 
with the projection$$\CC^2_{(p)} \stackrel{\pi_1}{\to} \CC_{(p)},$$where 
$ \CC^2_{(p)} \cong \F \otimes \Oc_p$, $\CC_{(p)} \cong T^1_{C,p}$ and $\pi_1 ((x,y)) = x$ (for more 
details, cf. \S 3 in \cite{F2}).
}

\end{remark}

\section{Association of elements 
in $T_{[s]}{\V}_{\delta}({\F})$ 
to singular surfaces in $| \Ii_{C/X} \otimes c_1(\F) |$ on $X$}\label{S:c}

The connection between singular schemes defined in Sections \ref{S:a} and \ref{S:b} can be 
further analyzed in order to determine interesting geometric interpretations of first-order 
deformations given by sections in $H^0(X, {\F}^{\Sigma})$. 

With notation as in \S \ref{S:a} and in Assumption \ref{ass:main}, 
let $[s] \in {\V}_{\delta}({\F})$, $C = V(s)$, $\Sigma = Sing(C)$. 
Denote by $\La:= c_1(\F) \in Pic(X)$. 

By \eqref{eq:koszul}, one has: 
\begin{equation}\label{eq:tspaces}
T_{[s]} ({\V}_{\delta}({\F})) \subset T_{[s]} ({\Pp}(H^0(X, {\F})) ) \cong \frac{H^0(X, {\F})}{H^0(X, \Oc_X)} 
\hookrightarrow H^0 (X, \Ii_{C/X} \otimes \La).
\end{equation}Therefore, first-order deformations of $[s]$ in ${\Pp}(H^0(X, {\F}))$, as well as in 
the Severi variety ${\V}_{\delta}({\F})$, 
can be related to suitable divisors moving in the linear system $|\La|$ on 
$X$ and containing the nodal curve $C$.

Indeed, we have: 

\begin{proposition}\label{prop:corr2}
Let $X $ be a smooth projective threefold. Let $\F$ be a globally generated rank-two vector
bundle on $X$ and let $\La = c_1(\F)$. Let $\delta$ be a positive integer, 
$[s] \in {\V}_{\delta}({\F})$ and $C=V(s)$ be the corresponding 
irreducible, nodal curve in $X$. Denote by $\Sigma$ the set of nodes of $C$. 
Let $\F^{\Sigma}$ be the sheaf on $X$ defined in Proposition \ref{prop:corr1}. 

Then: 

\begin{itemize}
\item[(i)] Each global section $s' \in H^0(X, {\F}^{\Sigma}) \setminus < s >$ determines 
a divisor $ V( s \wedge s') \in \; | \Ii_{C/X} \otimes \La| $ in X which 
is singular along $\Sigma$. 

\noindent
Precisely, 
\begin{equation}\label{eq:propcorr2}
s' \in  H^0(X, {\F}^{\Sigma})  \Leftrightarrow V (s \wedge s') \; {\rm contains} \; C \; {\rm 
and \; it \; is \; singular \; along} \; \Sigma.
\end{equation}

\item[(ii)] The singularities of $V(s \wedge s')$ are along $\Sigma$ and along the (possibly empty) 
intersection scheme  $C \cap V(s')$. 
\end{itemize}
\end{proposition}
\begin{proof} 
(i) Consider $s' \in H^0(X, {\F}^{\Sigma})$. By Proposition \ref{prop:corr1} and 
Remark \ref{rem:localdescr}, this corresponds to a global section 
of $\F$ which is in the kernel of the map $\mu_X$ - i.e. a global section 
which goes to zero in the composition $\F \to\!\!\to \F|_C \to\!\!\to T^1_C$ in 
diagram \eqref{eq:(*1)}. 

Since the situation is local, we may work locally around each node, in some 
open subset where $\F$ trivializes. Thus, fix 
a node $p \in \Sigma$ and a suitable neighborhood $U = U_p$ of $p$, whose local 
coordinates are denoted by $(x_1, x_2, x_3)$. 
We assume that $C$ is defined in $U$ by two equations $ f_1 = f_2 = 0$, 
where $s|_U = (f_1, f_2)$, $f_1, f_2 \in \Oc_X(U)$. 

Thus, the kernel of $\mu_X$ is given by those sections which are, at each node $p$, 
in the image of the Jacobian map
\begin{equation}\label{eq:jacobian}
{\T}_{{\CC}^3}|_C \stackrel{J(s)}{\longrightarrow} {\N}_{C/{\CC}^3} \to T^1_C
\end{equation}given 
by  
\begin{equation}\label{eq:jacobian2}
J(s) := \left( \begin{array}{ccc}
			\frac{\partial f_1}{\partial x_1} & \frac{\partial f_1}{\partial x_2} & 
                         \frac{\partial f_1}{\partial x_3}\\
                         \frac{\partial f_2}{\partial x_1} & \frac{\partial f_2}{\partial x_2} & 
                         \frac{\partial f_2}{\partial x_3}
		 	\end{array}
			\right) 
\end{equation}in $U$. 

Assume $s' \in H^0(X, {\F}^{\Sigma})$ and suppose $s'|_U = (g_1, g_2)$. 
Therefore, for each $p \in \Sigma$, $s' \in H^0(X, {\F}^{\Sigma})$ if and only if 
$s'(p) = (g_1(p), g_2(p)) \in \F_p$ is linear dependent on each pair 
$(\frac{\partial f_1}{\partial x_i}(p), \frac{\partial f_2}{\partial x_i}(p))$, $1 \leq i \leq 3$. 
This is equivalent to the following conditions:
$$
det \left( \begin{array}{ll}
	g_1(p) & g_2(p) \\
	\frac{\partial f_1}{\partial x_i}(p) & \frac{\partial f_2}{\partial x_i}(p) 
\end{array}	   
\right)= 0, \;{\rm for \; each}\; 1 \leq i \leq 3, 
$$i.e. 
\begin{equation}\label{eq:fico}
g_1 (p) \frac{\partial f_2}{\partial x_i}(p) - g_2 (p) \frac{\partial f_1}{\partial x_i}(p) = 0, \;{\rm for \; each}\; 
1 \leq i \leq 3.
\end{equation}

On the other hand observe that, since in particular 
$s' \in H^0(X, \F)$, then it defines a divisor 
in $| \La |$ on $X$ containing $C$. Indeed, consider 
$$\tau := (s, s') : \Oc_X \oplus \Oc_X \to \F,$$where
$$\tau :=\left( \begin{array}{ll}
	f_1 & g_1 \\
	f_2 & g_2 
\end{array}	   
\right)$$in the given open subset $U$. 
The degeneration locus of the map $\tau$ is given by$$V (det (\tau)) = V ( s \wedge s'),$$whose 
local equation in $U$ is given by
\begin{equation}\label{eq:loceq}
f_1g_2 - f_2 g_1 = 0.
\end{equation}Thus, since $\La = c_1(\F)$, 
$V(s \wedge s')$ corresponds to a divisor in $|\La|$ containing $C = V(s)$, i.e. it 
belongs to the linear system $| \Ii_{C/X} \otimes \La|$ on $X$. 

By the local equation of $V(s \wedge s')$ in $U$ and by the fact that $p \in C = V(s)$, 
we have that 
\begin{equation}\label{eq:contiloc}
\frac{\partial}{\partial x_i}(f_1 g_2 - g_1 f_2) (p) =  
g_1 (p) \frac{\partial f_2}{\partial x_i}(p) - g_2 (p) \frac{\partial f_1}{\partial x_i}(p), \; 
{\rm for \; each} \; 1 \leq i \leq 3. 
\end{equation}Therefore, $s' \in H^0(X, {\F}^{\Sigma})$ if and only if 
the associated divisor in $|\Ii_{C/X} \otimes \La|$ is singular at the nodes of $C$ 
(cf. \cite{BC}, for the case $X = \Pt$).

\vskip 10pt

\noindent
(ii) Let $q \in C$ be a point and let $s = (f_1, f_2)$, $s' = (g_1, g_2)$ be the 
local expressions of $s$ and $s'$ around $q$. Fix $(x_1, x_2, x_3)$ local coordinates 
of $X$ around $q$. Since $C = V(s)$ and $q \in C$, 
we have that \eqref{eq:contiloc} holds, for each $1 \leq i \leq \delta$. 

Assume that $q \in C \setminus \Sigma$ and that 
$s' (q) \neq (0,0)$; in this case, if \eqref{eq:contiloc} 
is equal to $0$ at $q$, for each $1 \leq i \leq 3$, we would have that 
$(g_1(q), g_2(q))$ is linear dependent on each pair 
\begin{equation}\label{eq:boh3}
(\frac{\partial f_1}{\partial x_i}(q), \frac{\partial f_2}{\partial x_i}(q)), \;  1 \leq i \leq 3.
\end{equation}In particular, the three pairs in \eqref{eq:boh3} would be 
linearly dependent. This is a contradiction; 
indeed, since $q \in C \setminus \Sigma$, the Jacobian map in \eqref{eq:jacobian} and \eqref{eq:jacobian2} 
is surjective at such a point $q$, i.e. $T^1_{C,q} = 0$. On the other hand, since 
${\N}_{C/X, q} \cong \Oc_{C,q}^{\oplus 2}$, then we must have that two of the three pairs 
in \eqref{eq:boh3} are linearly independent. 

This implies that $V(s\wedge s')$ 
cannot be singular outside $\Sigma \cup (C \cap V(s'))$. Indeed, 
in the other cases - i.e. either $q \in \Sigma$ or $q \in (C \cap V(s'))$ or both - 
it is easy to observe that \eqref{eq:contiloc} always vanishes at $q$, so that 
$V(s \wedge s')$ is singular at each such a point.
\end{proof}

By using the "divisorial" approach introduced in 
Theorem \ref{prop:3.fundamental}, we shall give in Theorem \ref{thm:fundamental} 
other interpretations of the equivalence in Proposition \ref{prop:corr2}, which 
highlights the deep connection between these apparently distinct 
approaches.

\begin{remark}\label{rem:propcorr2}
\normalfont{
From \eqref{eq:contiloc}, we see that among the global sections in 
$H^0(X, {\F}^{\Sigma})$, there are those global sections $s^*$ 
such that $s^* (p_i) = (0,0)$, for $p_i \in \Sigma$, 
$1 \leq i \leq \delta$. 

By the very definition of $\F^{\Sigma}$, 
we find that the inclusion$$\Ii_{\Sigma/X} 
\otimes {\F} \subseteq {\F}^{\Sigma}$$is not an isomorphism of sheaves on $X$. 
Indeed, the global sections of $H^0(X, {\F}^{\Sigma})$ are those satisfying a 
condition like \eqref{eq:fico} at each node in $\Sigma$; by \eqref{eq:contiloc}, 
condition \eqref{eq:fico} in particular holds if 
we consider global sections in the vector space $H^0(X , \Ii_{\Sigma/X} \otimes \F)$. 
Such a vector space has an expected codimension equal to $2 \delta$ in 
$H^0(X, \F)$ as it follows by the exact sequence:
\begin{equation}\label{eq:boh}
0 \to \Ii_{\Sigma/X} \otimes \F \to \F \to \F \otimes \Oc_{\Sigma} \cong \Oc_{\Sigma}^{\oplus 2} \to 0,
\end{equation}therefore, $\Ii_{\Sigma/X} \otimes \F$ is a proper subsheaf of ${\F}^{\Sigma}$. 

To sum up, surfaces in $X$ given by $V(s \wedge s')$, where $[s] \in {\V}_{\delta}({\F})$ and 
$s' \in H^0(X,\F^{\Sigma})$, are certainly 
singular along $\Sigma$ if the zero-locus $V(s')$ passes there; however, 
they can be also singular along $\Sigma$ even if $V(s')$ 
does not pass there. This happens since $C = V(s)$ is singular 
along $\Sigma$ and when \eqref{eq:fico} holds. 

There is also a "divisorial" interpretation of the fact that $H^0(X , \Ii_{\Sigma/X} \otimes \F)$ 
does not span the whole $H^0(X, {\F}^{\Sigma})$. From the 
correspondence of Theorem \ref{prop:3.fundamental} and from what stated in 
Corollary \ref{cor:propcorr1}, we know that the general section 
$s' \in  H^0(X, {\F}^{\Sigma})$ corresponds to a divisor $G_{s'}$ 
in the tautological linear system $ |\Oc_{\Pmc}(1) |$ 
on $\Pmc$, which simply passes through the scheme $\Sigma^1$ of $\delta$ 
rational double points of the divisor $G_s \in |\Oc_{\Pmc}(1) |$, where 
$G_s$ corresponds to the given section $[s] \in {\V}_{\delta}({\F})$ we started with. 
Recall that $G_s$ instead contains the whole $\pi$-fibres over 
$\Sigma$ and that 
it is singular along $\Sigma^1$. 

Thus, if one considers $ s^* \in H^0(X , \Ii_{\Sigma/X} \otimes \F)$, such an 
element corresponds to a divisor $G_{s^*} \in | \Ii_{\Lambda/\Pmc} \otimes \Oc_{\Pmc}(1) |$, 
where $\Lambda = \bigcup_{p_i \in \Sigma} L_{p_i}$. In this case, we have 
\begin{equation}\label{eq:boh1} 
0 \to \Ii_{\Lambda/\Pmc} \otimes \Oc_{\Pmc}(1) \to \Oc_{\Pmc}(1) \to 
\Oc_{\Pmc}(1)\otimes \Oc_{\Lambda} \cong \bigoplus_{i=1}^{\delta} \Oc_{L_{p_i}}(1) \to 0,
\end{equation}where $\Sigma= \{p_1, \ldots, p_{\delta} \}$. Observe that 
$$\Oc_{L_{p_i}}(1) \cong \Oc_{\Pp^1}(1), \; \forall \; 1 \leq i \leq \delta;$$obviously, 
$| \Ii_{\Lambda/\Pmc} \otimes \Oc_{\Pmc}(1) |$ is properly contained in 
$ |\Ii_{\Sigma^1/\Pmc} \otimes \Oc_{\Pmc}(1) |$, since the general element of the 
latter linear system 
simply passes through $\Sigma^1$ but does not contain the whole scheme $\Lambda$ as each element 
of $|\Ii_{\Lambda/\Pmc} \otimes \Oc_{\Pmc}(1) |$ does. 
Therefore, $ |\Ii_{\Sigma^1/\Pmc} \otimes \Oc_{\Pmc}(1) |$ 
has expected codimension equal to $2 \delta$ 
in $|\Oc_{\Pmc}(1) |$ as we found by using \eqref{eq:boh}.

For completeness sake, we conclude by observing 
that, in this correspondence, the subsheaf$$\Ii_{C/X} \otimes {\F} \subset \F^{\Sigma}$$gives global 
sections which are related to divisors on $\Pmc$ belonging to 
$| \Ii_{\mathbb{F}/\Pmc} \otimes \Oc_{\Pmc}(1) |$, where $\mathbb{F} = {\Pp}^1_C$ 
is the ruled surface contained in $\Pmc$ with $\pi$-fibres over the 
base curve $C$.
}
\end{remark}

\section{Connection among various singular subschemes of $X$ and of $\Pmc$}\label{S:d} 

The aim of this section is to study in details the deep connection among: 
\begin{itemize}
\item a nodal section $[s] \in {\V}_{\delta}({\F})$ on $X$,
\item the corresponding singular divisor $G_s$ in $|\Oc_{\Pmc}(1) |$ on $\Pmc$ having $\delta$ rational double points over 
the nodes of $C = V(s)$, 
\item for any $s' \in H^0 (X, {\F}^{\Sigma}) \setminus <s>$, 
the singular surface $V(s \wedge s') \subset X$ belonging to the linear system 
$| \Ii_{C/X} \otimes \La|$ on $X$, 
\item the singular "complete intersection" surface in $\Pmc$, 
${\Ss}_{s,s'} := G_s \cap G_{s'}$, dominating $V(s \wedge s')$ (cf. Remark \ref{rem:2propcorr2}). 
\end{itemize}

As a consequence of the analisys given 
in Sections \ref{S:b}, \ref{S:c}, we have the following:

\begin{theorem}\label{thm:fundamental}
Let $X $ be a smooth projective threefold. Let $\F$ be a globally generated rank-two vector
bundle on $X$ and let $\delta$ be a positive integer. Fix $[s] \in {\V}_{\delta}({\F})$ and
let $C=V(s)$ be the corresponding nodal curve on $X$. 
Denote by $\Sigma$ the set of nodes of $C$. 

Let ${\mathcal P} := {\Pp}_{X}({\F})$ be the
projective space bundle, ${\Oc}_{\mathcal P}(1)$ be its tautological line bundle 
and $\Sigma^1 \subset \Pmc$ be the zero-dimensional scheme of $\delta$-rational 
double points of the divisor $G_s \in |{\Oc}_{\mathcal P}(1)|$ which 
corresponds to $s \in H^0(X, \F)$ (cf. Theorem \ref{prop:3.fundamental}). 
Let $\F^{\Sigma}$ be the subsheaf of $\F$ defined 
as in Proposition \ref{prop:corr1}.

Then, the following conditions are equivalent:
\begin{itemize}
\item[(i)] $s' \in H^0(X, {\F}^{\Sigma})\setminus <s>$;
\item[(ii)] $V (s \wedge s') \subset X$ is a surface which 
contains $C$ and which is singular along $\Sigma$;
\item[(iii)] the divisor $G_{s'}$ passes through $\Sigma^1$
\item[(iv)] the surface $\Ss_{s,s'}:= G_s \cap G_{s'} \subset \Pmc$ is singular 
along $\Sigma^1$; 
\end{itemize} 
\end{theorem}
\begin{proof}Some of the implications are already proved in the previous results. Indeed:

\vskip 5 pt

\noindent
$(i) \Leftrightarrow (ii)$: this has already been 
proved in Proposition \ref{prop:corr2}.

\vskip 5 pt

\noindent
$(ii) \Leftrightarrow (iii)$: By the very 
definition of $\F^{\Sigma}$, it is obvious that 
$(iii)$ is equivalent to $(i)$; therefore, from the step above we would have finished. 
Anyhow, we give a direct proof of the equivalence of this two conditions, since 
it highlights the strict connection between the two apparently different approaches and it allows 
to give another interpretation of the equivalence in \eqref{eq:propcorr2}. 

To this aim, let 
$[s] \in {\V}_{\delta}({\F})$ and let $G_s$ be the corresponding 
divisor in $\Pmc$ which is singular along $\Sigma^1$. As usual, 
take $p \in \Sigma$, $U = U_p$ an open subset such that 
$U \cap (\Sigma \setminus \{ p \}) = \emptyset$, 
where $\F$ trivializes and whose local coordinates 
are $\underline{x} = (x_1, x_2, x_3)$, such that 
$\underline{x}(p) = \underline{0} \in U \cong 
\mathbb{A}^3$. If $s|_U = (f_1,f_2)$, 
we recall that the local equation of $G_s$ in $\pi^{-1} (U)$ is 
$$F(x_1, x_2, x_3, u ,v):= u f_1 + v f_2$$ Therefore, we have 
a singular point on $G_s$ in $\pi^{-1}(U)$ if, and only if, there exists a solution of
\begin{equation}\label{eq:system1}
F = \frac{\partial F}{\partial x_1} = \frac{\partial F}{\partial x_2} = 
\frac{\partial F}{\partial x_3} = 
\frac{\partial F}{\partial u} = \frac{\partial F}{\partial v} = 0. 
\end{equation}Observe that \eqref{eq:system1} is equivalent to
\begin{equation}\label{eq:system2}
\begin{array}{l}
uf_1+vf_2 = \frac{\partial f_1}{\partial x_1} u + \frac{\partial f_2}{\partial x_1} v = 
\frac{\partial f_1}{\partial x_2} u + \frac{\partial f_2}{\partial x_2} v = \\
= \frac{\partial f_1}{\partial x_3} u + \frac{\partial f_2}{\partial x_3} v = 
f_1 = f_2 = 0.
\end{array}
\end{equation}By the last two equations, we find as in Theorem \ref{prop:3.fundamental} that 
the singular point of $G_s$ must be on the $\pi$-fibre 
over a point of $C$. Let $q \in U \cap C$ be such a point and let $L = L_q$ be 
this fibre. 

We can restrict the system \eqref{eq:system2} to $L$. We thus get:
\begin{equation}\label{eq:system3}
\frac{\partial f_1}{\partial x_1}(q) u + \frac{\partial f_2}{\partial x_1}(q) v = 
\frac{\partial f_1}{\partial x_2}(q) u + \frac{\partial f_2}{\partial x_2}(q) v = 
\frac{\partial f_1}{\partial x_3}(q) u + \frac{\partial f_2}{\partial x_3}(q)v = 0.
\end{equation}Therefore, there exists a solution $[u,v] \in L \cong \Pp^1$ 
if, and only if, \eqref{eq:system3} has rank less than or equal to one. This 
is equivalent to saying that
\begin{equation}\label{eq:propvectors}
(\frac{\partial f_1}{\partial x_1}(q) , \frac{\partial f_1}{\partial x_2}(q), 
\frac{\partial f_1}{\partial x_3}(q)) = \lambda (\frac{\partial f_2}{\partial x_1}(q) , 
\frac{\partial f_2}{\partial x_2}(q), 
\frac{\partial f_2}{\partial x_3}(q)),
\end{equation}for some $\lambda \in \CC^*$; this is equivalent to:
\begin{equation}\label{eq:system4}
\begin{array}{l}
\frac{\partial f_1}{\partial x_1}(q) \frac{\partial f_2}{\partial x_2}(q) - 
\frac{\partial f_2}{\partial x_1}(q) \frac{\partial f_1}{\partial x_2}(q) = \\ 
\frac{\partial f_1}{\partial x_1}(q) \frac{\partial f_2}{\partial x_3}(q) - 
\frac{\partial f_2}{\partial x_2}(q) \frac{\partial f_1}{\partial x_3}(q) = \\
\frac{\partial f_1}{\partial x_2}(q) \frac{\partial f_2}{\partial x_3}(q) - 
\frac{\partial f_1}{\partial x_3}(q) \frac{\partial f_1}{\partial x_2}(q) = 0.
\end{array}
\end{equation}In such a case, we have:
\begin{equation}\label{eq:solution}
\begin{array}{rcl}
[u,v]  & = & [-\frac{\partial f_2}{\partial x_1}(q), \frac{\partial f_2}{\partial x_1}(q)] 
= [-\frac{\partial f_2}{\partial x_2}(q), \frac{\partial f_1}{\partial x_2}(q)] \\ 
& = & [-\frac{\partial f_2}{\partial x_3}(q), \frac{\partial f_2}{\partial x_3}(q)]  
= [- \lambda, 1].
\end{array}
\end{equation}Therefore, $G_s$ is singular at $[- \lambda, 1] \in L$ if, and only if, 
$C = V(s) $ has a node at $q \in U$. This implies that $q=p$, since $p$ was the only node 
in $U$ by assumption; thus $L = L_p$ and - once again - 
the singularities of $G_s$ are on the $\pi$-fibres of the nodes of $C$.

Let now $s' \in H^0(X, \F^{\Sigma})$. Assume that 
$s'|_U = (g_1, g_2)$, for some $g_1, g_2 \in \Oc_{X}(U)$. 
Then, $G_{s'}$ passes through the singular point of $G_s$ along 
$L_p$, i.e. $[-\lambda, 1]$ if, and only if, 
\begin{equation}\label{eq:propvectors2}
[-g_2(p), g_1(p)] = [-\lambda, 1].
\end{equation}This means that $[-g_2(p), g_1(p)]$ is a solution of \eqref{eq:system3} (where 
$q=p$, as proved above). This is exactly equivalent to \eqref{eq:fico} and so to 
the fact that the surface $V( s \wedge s')$ is singular at $p$.

\vskip 5 pt

\noindent
$(iii) \Leftrightarrow (iv)$: trivial consequence of the fact that $G_s$ is always singular at 
$\Sigma^1$.

\end{proof}

\begin{remark}\label{rem:2propcorr2} 
\normalfont{
The proof of the equivalence of conditions $(ii)$ and 
$(iii)$ in Theorem \ref{thm:fundamental} shows that 
the local computations on $X$ in Proposition \ref{prop:corr2} are exactly equivalent 
to the local computations on $\Pmc$ via the correspondence introduced in 
Theorem \ref{prop:3.fundamental} and Proposition \ref{prop:corr1}.

In particular, from \eqref{eq:solution} and \eqref{eq:propvectors2}, it follows that 
\begin{equation}\label{eq:propvectors3}
[-g_2(p), g_1(p)] = [-\frac{\partial f_2}{\partial x_i}(q), 
\frac{\partial f_2}{\partial x_i}(q)], 
\; {\rm for \; each} \; 1 \leq i \leq 3. 
\end{equation}This means that $(g_1(p), g_2(p))$ linearly depends on 
each pair 
$(\frac{\partial f_2}{\partial x_i}(q), \frac{\partial f_2}{\partial x_i}(q))$, 
$1 \leq i \leq 3$, as already obtained 
in \eqref{eq:boh3}. 

It is interesting to give also a direct proof of the equivalence of conditions 
$(ii)$ and $(iv)$, in order to relate it with what 
observed in Remark \ref{rem:propcorr2}. Thus:

\vskip 5pt

\noindent
$(ii) \Rightarrow (iv)$: Since the computations are local, we can fix 
one node $p \in \Sigma$. Since $V (s \wedge s')$ is singular along $\Sigma$ 
by assumption, from \eqref{eq:contiloc} it follows that either 
$s'$ passes through $p$ or $s'(p)$ is proportional to each pair as in \eqref{eq:boh3} evaluated 
at $p$. In the former case, we have that $G_{s'} \in |\Ii_{L_p/\Pmc} \otimes \Oc_{\Pmc}(1) |$, where 
$L_p$ is the $\pi$-fibre over $p \in \Sigma$; in the latter case, by the very definition 
of $\F^{\Sigma}$ and by \eqref{eq:jacobian} and \eqref{eq:jacobian2}, we have that 
$G_{s'} \in |\Ii_{p^1} \otimes \Oc_{\Pmc}(1) |$, where $p^1 \in \Sigma^1$ is the corresponding 
point to $p \in \Sigma$. In any case, $G_{s'}$ passes through $p^1$. 

If we globalize this approach, in any case, $G_{s'}$ passes through $\Sigma^1$. Now, 
since $\Sigma^1 \subseteq Sing (G_s)$, then it obviously follows that 
$\Ss_{s,s'}:= G_s \cap G_{s'}$ is singular along $\Sigma^1$.

\vskip 5 pt

\noindent
$(iv) \Rightarrow (ii)$: Once again, we can focus on one 
of the nodes in $\Sigma$, say $p$. Take $U = U_p$ an open neighborhood of $p$ in 
$X$ where $\F$ trivializes. Then, we can assume that the local expression of 
$s$ in $U$ is $s|_U = (f_1, f_2),$ where $f_1, f_2 \in \Oc_X(U)$. Denote by $\mop$ the 
maximal ideal of the point $p$ in the stalk $\Oc_{X,p}$. 
Since by assumption $[s] \in {\V}_{\delta}({\F})$ and $p \in \Sigma$, 
we can assume that the reduction of $s$ in $\F \otimes (\mop/\mop^2)$ is 
$(1,0)$. This means that if we consider homogeneous coordinates $[u,v]$ on the 
$\pi$-fibre $L_p \cong \Pp^1$ over p, the corresponding 
rational double point $p^1$ for $G_s$ on $L_p$ 
has coordinates $[0,1]$ on such a line.

Similarly, take $s'|_U = (g_1, g_2),$ where $g_1, g_2 \in \Oc_X(U)$. Since by assumption 
$\Ss_{s,s'} = G_s \cap G_{s'}$ is singular at $p^1 = [0,1]$, in particular $G_{s'}$ passes 
through $p^1$. Therefore, we can assume 
that the reduction of $s'$ in $\F \otimes (\Oc_{X,p}/\mop)$ is $(a,0)$. If $a = 0$, 
this means that $G_{s'}$ contains $L_p$; otherwise, as in \eqref{eq:loccomp0}, the local equation 
of $G_{s'}$ is given by $\{a u = 0\}$, so that the intersection point between $G_{s'}$ and 
$L_p$ is indeed $p^1= [0,1]$. 

In any case, we have that $g_2 \in \mop$ and 
$g_1 = a + j_1$, where $j_1 \in \mop$. Analogously, we have that $f_1 \in \mop$ and 
$f_2 \in \mop^2$. Therefore, $$det \left( \begin{array}{ll}
	f_1 & f_2 \\
	a+j_1 & g_2 
\end{array}	   
\right) = f_1 g_2 - f_2 (a+ j_1) \in \mop^2.$$This implies that $V (s \wedge s')$ is singular 
at $p \in \Sigma$.

Recall that, in Remark \ref{rem:propcorr2} we observed that 
surfaces in $X$ given by $V(s \wedge s')$, with $s' \in H^0(X, {\F}^{\Sigma})$ 
are certainly singular along $\Sigma$ if the zero-locus $V(s')$ passes there; however, 
they can be also singular along $\Sigma$ even if $V(s')$ 
does not passes there, precisely when \eqref{eq:fico} holds at each point of $\Sigma$. 
From the correspondence between $V(s \wedge s')$ and $\Ss_{s,s'}$ 
we see that, in the former case, the surface $\Ss_{s,s'}$ has to contain 
$\Lambda = \pi^{-1}(\Sigma) = \bigcup_{i=1}^{\delta}L_{p_i}$, whereas in the latter, $\Ss_{s,s'}$ has to pass 
through the point $p_i^1 \in L_{p_i}$, for each $p_i \in \Sigma$, which is singular 
for $G_s$ so - a fortiori - for $\Ss_{s,s'}$. 

In any case, differently from $V(s \wedge s')$, the surface 
$\Ss_{s,s'}$ always contains $\Sigma^1$ if $s' \in H^0(X, {\F}^{\Sigma})$. 
Notice also that $\Ss_{s,s'}$ is a "complete intersection" in $\Pmc$ which dominates 
$V(s \wedge s')$. Indeed, if $s' \in H^0(X, {\F}^{\Sigma})$ is such that either 
$V(s') \cap \Sigma \neq \emptyset$ or $V(s') \cap C \neq \emptyset$ then, as above, 
the $\pi$-fibres over these intersection schemes 
are properly contained in $\Ss_{s,s'}$. In particular, 
if $s'$ is the general section in $H^0(X, {\F}^{\Sigma})$ such that $V(s') \cap C = \emptyset$ 
then, by the Zariski Main Theorem, $\Ss_{s,s'}$ is isomorphic via $\pi$ to the surface 
$V(s \wedge s')$. 
}
\end{remark}

Observe that by Theorem \ref{thm:fundamental}, we can improve 
the statement of Corollary \ref{cor:propcorr1}. 

\begin{coroll}\label{cor:explanetion} $s + \epsilon s' \in T_{[s]} ({\V}_{\delta}({\F}))$, s.t. $\epsilon^2 = 0$ 
$\Leftrightarrow$ $V(s \wedge s')$ is singular along $\Sigma$ $\Leftrightarrow $
$\G_{s'} \in |\Ii_{\Sigma^1/\Pmc} \otimes \Oc_{\Pmc}(1) | $ $\Leftrightarrow $ 
$\Ss_{s,s'}:=\G_s \cap \G_{s'} \subset \Pmc$ is singular along $\Sigma^1$.
\end{coroll}

\begin{remark}\label{rem:explanetion}
\normalfont{
By using Theorem \ref{thm:fundamental}, we can also 
give another interpretation 
of the equivalence of regularity conditions in Corollary \ref{rem:reg}. Recall that 
the map $\rho_{\Pmc}$ in \eqref{eq:regbis} is a standard 
restriction map. 

Therefore, $| \Oc_{\Pmc}(1) |$ does not separate $\Sigma^1$ if, and only if, each 
divisor in $| \Oc_{\Pmc}(1) |$ passing through all but one point $p_j^1$ 
of $\Sigma^1$ passes also through the point $p_j^1$, for some $1 \leq j \leq \delta$. 
By Theorems \ref{prop:3.fundamental} and \ref{thm:fundamental} and 
by Proposition \ref{prop:corr1}, this happens if, and only if, for each 
$[s] \in {\V}_{\delta}({\F})$ and for each 
$ s' \in H^0 (X, {\F}^{\Sigma})$, 
the surface $\Ss_{s,s'}= G_s \cap G_{s'}$ which is singular along all but one point $p_j^1$ 
of $\Sigma^1$ is singular also at the remaining 
point $p_j^1$, for some $1 \leq j \leq \delta$. This happens if, and only if, for each 
$[s] \in {\V}_{\delta}({\F})$ and for each 
$ s' \in H^0 (X, {\F}^{\Sigma})$, 
the surface $V(s \wedge s') \subset X$ which is singular along all but one point $p_j$ 
of $\Sigma$ is singular also at the remaining 
point $p_j$, for some $1 \leq j \leq \delta$; this is equivalent to the non-surjectivity 
of the map $\mu_X$, since the section $s' \in H^0(X, \F)$ which vanishes in the 
composition $\F \to \!\! \to \F|_C \to \!\! \to \Oc_{\Sigma \setminus \{ p_j\}}$ also vanishes 
in the composition $\F \to \!\! \to \F|_C \to \!\! \to \Oc_{\{ p_j\}}$. 
}
\end{remark}

To sum up, the regularity condition for points of the scheme ${\V}_{\delta}({\F})$ on $X$ 
not only translates into the independence of conditions imposed by a $0$-dimensional scheme 
on the tautological linear system $|\Oc_{\Pmc}(1) |$ on $\Pmc$, but also to the 
independence of {\em singularity conditions} imposed on suitable families of singular surfaces 
in $X$ as well as in $\Pmc$.

\section{$\Pmc_{\delta}$-Severi varieties of singular 
divisors on $\Pmc$}\label{S:ef}

The aim of this section is to study 
families of singular divisors in $\Pmc$, which are 
related to Severi varieties ${\V}_{\delta}({\F}) $ of nodal sections 
of $\F$ on $X$.

\begin{definition}\label{def:rdelta}
Given $\Pmc$ and $\delta$ as in Assumption \ref{ass:main}, Theorem 
\ref{prop:3.fundamental} and Proposition \ref{prop:corr1}, 
consider the scheme 
\begin{equation}\label{eq:rdelta}
{\mathcal R}_{\delta} (\Oc_{\Pmc}(1)) 
:=  \{G_s \in |\Oc_{\Pmc}(1) | \; {\rm s.t.} \; [s] \in {\V}_{\delta}({\F}) \}. 
\end{equation}These schemes parametrize families of 
divisors in the tautological linear system of $\Pmc$, $|\Oc_{\Pmc}(1) |$, 
which are irreducible and have only $\delta$ rational double points as the only 
singularities. For brevity sake, these 
will be called $\Pmc_{\delta}$-{\em Severi varieties}.
\end{definition}

By the above definition and by Assumption \ref{ass:main}, from now on we 
consider $\delta \leq \; min\{ h^0(\F) -1, \frac{1}{2}deg (\La \otimes \omega_X \otimes \Oc_D) + 1\}$ 
and ${\mathcal R}_{\delta} (\Oc_{\Pmc}(1)) \neq \emptyset$. 

It is clear that:
\begin{equation}\label{eq:expdim}
expdim ({\mathcal R}_{\delta} (\Oc_{\Pmc}(1))) = dim (|\Oc_{\Pmc}(1)|) - \delta;
\end{equation}indeed, imposing a rational double point gives 
at most $5$ conditions on $|\Oc_{\Pmc}(1)|$; each such point 
varies on any of the $\pi$-fibre over $X$.

As in Definition \ref{def:0}, from \eqref{eq:expdim} it is natural to give the 
following:

\begin{definition}\label{def:expdim}
Let $[G_s] \in {\mathcal R}_{\delta} (\Oc_{\Pmc}(1))$.
Then $[G_s]$ is said to be a {\em regular point} of ${\mathcal R}_{\delta} (\Oc_{\Pmc}(1))$ if:
\begin{itemize}
\item[(i)] $[G_s] \in {\mathcal R}_{\delta} (\Oc_{\Pmc}(1))$ is a smooth point, and
\item[(ii)] $dim_{[G_s]}({\mathcal R}_{\delta} (\Oc_{\Pmc}(1))) = expdim ({\mathcal R}_{\delta} (\Oc_{\Pmc}(1)))=
 dim (|\Oc_{\Pmc}(1)|) - \delta$. 
\end{itemize}The $\Pmc_{\delta}$-Severi variety 
${\mathcal R}_{\delta} (\Oc_{\Pmc}(1))$ is said to be {\em regular} if it is 
regular at each point. 
\end{definition}

As in Theorem \ref{prop:3.fundamental}, in order 
to find regularity conditions for a given point $[G_s] \in {\mathcal R}_{\delta} (\Oc_{\Pmc}(1))$, 
it is crucial to give a description of the tangent space at 
$[G_s]$ to the given $\Pmc_{\delta}$-Severi variety.

\begin{theorem}\label{thm:tangentspace} Let 
$[G_s] \in {\mathcal R}_{\delta} (\Oc_{\Pmc}(1))$ on $\Pmc$ and let $\Sigma^1 $ be 
the zero-dimensional scheme of the $\delta$-rational double points of $G_s \subset \Pmc$.
Then, we have:
\begin{equation}\label{eq:tangentspace}
T_{[G_s]} ({\mathcal R}_{\delta} (\Oc_{\Pmc}(1))) \cong \frac{H^0(\Ii_{\Sigma^1/\Pmc} \otimes \Oc_{\Pmc}(1))}{<G_s>}.
\end{equation}

\noindent
In particular, if $\epsilon \in \CC[T]/(T^2)$, then:
$$G_s + \epsilon \; G_{r} \in T_{[G_s]} ({\mathcal R}_{\delta} (\Oc_{\Pmc}(1)) \Leftrightarrow  
G_r \in |\Oc_{\Pmc}(1)) | \; {\rm and} \;  \Sigma^1 \subset G_r.$$

\end{theorem}
\begin{proof}
The divisor $G_s \subset \Pmc$, related to the 
point $[G_s] \in {\mathcal R}_{\delta}(\Oc_{\Pmc}(1))$, corresponds to a section 
$[s] \in  {\V}_{\delta}({\F})$ on $X$. 

Therefore, as in the proof of Theorem \ref{prop:3.fundamental}, we may locally work 
around a node of $C = V(s) \subset X$.  Let $p \in \Sigma = Sing(C)$ be such a node and let $U = U_p \subset X$ be 
an affine open set containing $p$, where the vector bundle ${\F}$ trivializes. 
We thus can choose local coordinates 
$\underline{x} = (x_1, x_2, x_3)$ on $U \cong {\mathbb A}^3$ and homogeneous
coordinates $[u,v] \in {\Pp}^1$, such 
that $\underline{x}(p) = (0,0,0)$, $s|_{U} = (x_1x_2, \; x_3)$ and the 
local equation of $G_s$ in $\pi^{-1}(U) \cong  U_p \times {\Pp}^1$ 
is given by $u x_1 x_2 + v x_3 = 0$ (cf. \eqref{eq:loccomp0}). 

Recall that in the open chart where $v \neq 0$, $G_s$ is smooth, whereas, in the open chart where $u \neq 0$, 
we see that the local equation of $G_s $ in ${\mathbb A}^3 \times {\mathbb A}^1 \cong \mathbb{A}^4$ 
is 
\begin{equation}\label{eq:palle}
G_s = V( x_1 x_2 + x_3 t), \; {\rm where} \; t= \frac{v}{u}.
\end{equation}This is the equation of a quadric cone in $\mathbb{A}^4$ having 
vertex at the origin of the $\mathbb{A}^4$ having coordinates $(x_1, x_2, x_3, t)$. 

We can consider the Jacobian map of $G_s$ in this $\mathbb{A}^4$. This is given by:
\begin{displaymath}
\begin{array}{rcl}
\T_{\mathbb{A}^4 |_{G_s}} & \stackrel{J_{G_s}}{\longrightarrow} & {\N}_{G_s/\mathbb{A}^4} \\
\partial /\partial x_1 & \longrightarrow& x_2 \\
\partial / \partial x_2 & \longrightarrow & x_1 \\
\partial / \partial x_3 & \longrightarrow & t \\
\partial / \partial t  & \longrightarrow & x_3,
\end{array}
\end{displaymath}where ${\N}_{G_s/\mathbb{A}^4}$ is locally free of rank one 
on $G_s$. It is then clear that $J_{G_s}$ is surjective except at the origin $\underline{0} = (0,0,0,0)$. 

By the local computations we analytically get: 
$$coker (J_{G_s}) \cong \frac{\CC[[x_1, x_2, x_3, t]]/(x_1 x_2 + x_3 t)}{(x_1, x_2, x_3, t)} \cong \CC;$$

Globally speaking, given $G_s \subset \Pmc$, whose singular scheme is $\Sigma^1$, 
we have the exact sequence of sheaves on $G_s$:
\begin{equation}\label{eq:palle1}
\T_{\Pmc|_{G_s}} \stackrel{J_{G_s}}{\to} {\N}_{G_s/\Pmc} \to T^1_{G_s} \to 0,
\end{equation}where $T^1_{G_s} $ is a sky-scraper sheaf supported on $\Sigma^1$ and 
of rank one at each point. 

As in \eqref{eq:T1} for nodal curves, denote by $\N'_{G_s}$ the kernel of $J_{G_s}$ in 
\eqref{eq:palle1}. This is the so-called {\em equisingular sheaf}, whose global sections 
give equisingular first-order deformations of $G_s$ in $\Pmc$. 

By standard exact sequences, one sees that there is an injection
$$ \frac{H^0(\Pmc, \Ii_{\Sigma^1/\Pmc} \otimes \Oc_{\Pmc}(1))}{H^0(\Pmc, \Oc_{\Pmc})} \hookrightarrow 
H^0(G_s, \N'_{G_s}),$$which is an isomorphism when $\Pmc$ - equivalently $X$ - is regular, 
i.e. $h^1(\Pmc, \Oc_{\Pmc}) = 0$. Therefore, the vector space on the left-hand-side of the 
injection actually parametrizes 
equisingular first-order deformations of $G_s$ in $|\Oc_{\Pmc}(1)|$. 
\end{proof}

\begin{remark}\label{rem:nfi} 
\normalfont{Recall that when 
one studies classical Severi varieties of irreducible, $\delta$-nodal 
curves on a smooth projective surface, 
there is also a {\em parametric approach} for equisingular first-order 
deformations (cf., e.g \cite{CH} and \cite{S}). 

Precisely, let $S$ be an arbitrary smooth, projective surface, 
$|D|$ be a complete linear system on $S$, whose general element 
is assumed to be a smooth and irreducible curve, which is a divisor on $S$; one considers the Severi variety 
$V_{|D|, \delta} \subset |D|$, for any $0 \leq \delta \leq p_a(D)$, which parametrizes 
reduced, irreducible curves in $|D|$ having 
$\delta$-nodes as the only singularities. If $[C] \in V_{|D|, \delta}$, this point 
corresponds to a curve $C \sim D$ on S, 
such that $N:= Sing(C) \subset S$ is the $0$-dimensional scheme of its 
$\delta$ nodes and one can consider:
\begin{equation}\label{eq:palle2}
\begin{array}{clccl}
\tilde{C} &  & \subset & \tilde{S} & \\
\downarrow & \!\!\!\!^{\varphi_N} & & \downarrow & \!\!\!^{\mu_N} \\
C & & \subset & S & , 
\end{array}
\end{equation}where
\begin{itemize}
\item $\mu_N$ is the blow-up of $S$ along $N$, 
\item $\varphi_N$ is the normalization of $X$,
\item $\Ct$ is a smooth curve of (geometric) genus $g= g(\Ct)= p_a(D) - \delta$.
\end{itemize}

It is a standard result that $T_{[X]}(V_{| D |, \; \delta}) \cong \frac{H^0(S, \; \Ii_{N/S}(D))}{<C>}$ is 
isomorphic to a (proper) subspace of $H^0({\N}_{\varphi_N})$, where 
${\N}_{\varphi_N}$ is the {\em normal bundle to map} $\varphi_N$, which is the line bundle 
on the smooth curve $\Ct$ defined by: 
$$0 \to {\T}_{\Ct} \to \varphi^*({\T}_S) \to {\N}_{\varphi_N} \to 0. $$It is well-known 
that $H^0(\Ct, {\N}_{\varphi_N})$ parametrizes equisingular first-order deformations of $C$ in $S$ and that 
the subspace mentioned above coincides with the whole space when $S$ is a regular surface. 

Therefore, for irreducible nodal curves on surfaces, the parametric approach coincides with 
the {\em Cartesian approach}, which makes use of the equisingular sheaf $\N'_C$ defined 
by$$ 0 \to \N'_C \to {\N}_{X/S} \to T^1_{X} \to 0.$$Indeed, in this case (and only in this case) 
one has $\N'_C \cong \varphi_*({\N}_{\varphi_N})$. 

Even if the elements of ${\mathcal R}_{\delta}(\Oc_{\Pmc}(1))$ are divisors 
in $\Pmc$ - as curves on surfaces - the same does not occur for these families. For simplicity, 
assume that $X$ - and so $\Pmc$ - is a regular threefold, i.e. $h^1(X, \Oc_X) = 0$. Let 
$$\mu_{\Sigma^1} : \tilde{\Pmc} \to \Pmc$$be the blow-up of $\Pmc$ along 
$\Sigma^1$ and$$\varphi_{\Sigma^1} : \tilde{G}_s \to G_s$$the desingularization of $G_s$, which 
is induced by $\mu_{\Sigma^1}$, by a diagram similar to the one in \eqref{eq:palle2} and by the fact 
that $\Sigma^1$ is a scheme of ordinary double points for $G_s$. Let $B := 
\Sigma_{i=1}^{\delta} E_i$ be the 
$\mu_{\Sigma^1}$-exceptional divisor. Thus, 
$$\mu_{\Sigma^1}^* (G_s) =  \tilde{G}_s + 2 B, \;\; \mu_{\Sigma^1}^* (K_{\Pmc}) = K_{\tilde{\Pmc}} - 3 B.$$By the 
exact sequence:$$0 \to {\T}_{\tilde{G_s}} \to \varphi_{\Sigma^1}^*({\T}_{\Pmc}) 
\to {\N}_{\varphi_{\Sigma^1}} \to 0$$and 
by the adjunction formula on $\tilde{\Pmc}$, we get that: 
\begin{equation}\label{eq:palle4}
{\N}_{\varphi_{\Sigma^1}} \cong \Oc_{\tilde{G_s}} (\mu_{\Sigma^1}^* (G_s) + B).
\end{equation}Tensoring by $\Oc_{\tilde{G_s}} (\mu_{\Sigma^1}^* (G_s) + B) $ the exact 
sequence$$ 0 \to \Oc_{\tilde{\Pmc}} (- \mu_{\Sigma^1}^* (G_s) + 2 B) \to 
 \Oc_{\tilde{\Pmc}} \to  \Oc_{\tilde{G_s}} \to 0,$$by the regularity of $\Pmc$ and 
by Fujita's Lemma, we see that $H^0({\N}_{\varphi_{\Sigma^1}})$ is not isomorphic to $H^0(\N'_{G_s})$, i.e. 
the first-order deformations given by general vectors in $H^0({\N}_{\varphi_{\Sigma^1}})$ 
are not equisingular.
}
\end{remark}

To conclude the section denote, as in \eqref{eq:regbis}, by$$\rho_{\Pmc} : H^0(\Pmc , \Oc_{\Pmc}(1)) 
\to H^0(\Pmc, \Oc_{\Sigma^1})$$the 
natural restriction map. Then, from Theorem \ref{thm:tangentspace}, it immediately follows:

\begin{corollary}\label{cor:tangentspace}
With assumptions and notation as in Theorem \ref{thm:tangentspace}, we have:
\begin{displaymath}
\begin{array}{lcl}
[G_s] \in {\mathcal R}_{\delta} (\Oc_{\Pmc}(1)) \; {\rm is \; a \; regular \; point} &  \Leftrightarrow & 
\rho_{\Pmc} \;  {\rm is \; surjective} \\
 & \Leftrightarrow & [s] \in {\V}_{\delta}({\F}) \; {\rm is \; a \; regular \; point}\\
 &  & {\rm (in \;  the \; sense \; of \; Definition \; \ref{def:0}).}
\end{array}
\end{displaymath}
\end{corollary}
\begin{proof} The first equivalence is a direct consequence of 
\eqref{eq:expdim} and Theorem \ref{thm:tangentspace}. The other 
follows from Corollary \ref{rem:reg}.
\end{proof}

\section{Some uniform regularity results for ${\V}_{\delta}({\F})$ and 
${\mathcal R}_{\delta} (\Oc_{\Pmc}(1)) $}\label{S:ghil}

In this section we first improve some regularity results of \cite{F2} 
for Severi varieties ${\V}_{\delta}({\F})$ of irreducible, 
$\delta$-nodal sections of $X$, then we use Corollary \ref{cor:tangentspace} 
to deduce regularity results also for $\Pmc_{\delta}$-Severi varieties 
${\mathcal R}_{\delta} (\Oc_{\Pmc}(1)) $ of irreducible divisors in $|\Oc_{\Pmc}(1) |$ on 
$\Pmc$. We find upper-bounds on the number $\delta$ of singular points 
which ensure the 
regularity of ${\V}_{\delta}({\F})$ as well as of ${\mathcal R}_{\delta} (\Oc_{\Pmc}(1)) $; these 
upper-bounds are shown to be almost sharp (cf. Remark \ref{rem:10}).

What we want to stress is the following fact: even if 
the regularity of the schemes ${\mathcal R}_{\delta} (\Oc_{\Pmc}(1)) $ on $\Pmc$ 
is defined by means of separation of suitable 
zero-dimensional schemes by the linear 
system $|\Oc_{\Pmc}(1)|$ on the fourfold $\Pmc$, one can avoid to consider this intricate 
situation. Indeed, it is well-known how difficult is to enstablish separation 
of points in projective varieties of dimension greater than or equal to three 
(cf. e.g. \cite{AS}, \cite{EL} and \cite{K}). 
In general some results can be found by using the technical tools 
of {\em multiplier ideals} together the Nadel and the Kawamata-Viehweg vanishing theorems (see, e.g. \cite{E}, 
for an overview). Anyhow, sometimes no answers are given by using these techniques. 

In our situation, thanks to the correspondence between 
${\V}_{\delta}({\F})$ on $X$ and 
${\mathcal R}_{\delta} (\Oc_{\Pmc}(1))$ on $\Pmc$, we deduce regularity 
conditions for the scheme ${\mathcal R}_{\delta} (\Oc_{\Pmc}(1))$ from those 
of the scheme ${\V}_{\delta}({\F})$.

From now on, let $X$ be a smooth projective threefold, 
$\E$ be a globally generated rank-two vector
bundle on $X$, $M$ be a very ample line bundle on $X$ and 
$k \geq 0$, $\delta >0 $ be integers.
With notation and assumptions as in Section \ref{S:a}, 
we shall always take$${\F} = {\E} \otimes M^{\otimes k}$$and consider the scheme
${\V}_{\delta}({\E} \otimes M^{\otimes k})$ on $X$.

By using Theorem \ref{prop:3.fundamental} and Corollary \ref{rem:reg}, here we determine
conditions on the vector bundle $\E$ and on the integer $k$ 
and uniform upper-bounds on the number of nodes
$\delta$ implying that each point of 
$ {\V}_{\delta}({\E} \otimes M^{\otimes k})$ is regular. Indeed, the next result is a 
generalization of Theorem 4.5 in \cite{F2}.

\begin{theorem}\label{thm:9bis}
Let $X$ be a smooth projective threefold, $\E$ be a globally generated rank-two vector bundle on 
$X$, $M$ be a very ample line bundle on $X$ and $k \geq 0$ and $\delta >0$ be integers. If 
\begin{equation}\label{eq:uniform}
\delta \leq k+1,
\end{equation}then ${\V}_{\delta}({\E} \otimes M^{\otimes k})$ is regular.
\end{theorem}
\begin{proof}
If $k =0$, then we consider $\delta = 1$; therefore, by the hypothesis on $\E$, it follows that 
$$H^0 (\E) \to H^0(\Oc_p^{\oplus 2}) $$is surjective, for each $p \in X$. By the description 
of the map $\mu_X$ in Remark \eqref{rem:localdescr}, this implies 
that ${\V}_1({\E})$ is regular at each point. 

When $k > 0$, first of all observe that, since $M$ is very ample, then 
$M^{\otimes k}$ separates any set 
$\Sigma$ of $\delta$ distinct point of $X$ with $\delta \leq k+1$. This is 
equivalent to saying that the restriction map
\begin{equation}\label{eq:prop7}
\rho_k: \; H^0(X, M^{\otimes k}) \to H^0({\Oc}_{\Sigma})
\end{equation}is surjective, for each such 
$\Sigma \subset X$. Thus, 
\eqref{eq:prop7} implies there exist global sections 
$\sigma_1, \ldots, \sigma_{\delta} \in H^0(X, M^{\otimes k})$ s. t. 
$${\sigma}_i(p_j) = \underline{0} \in {\CC}^{\delta}, \; {\rm if \; i \neq j }, \;{\rm and} \;  
\sigma_i(p_i) = (0, \ldots, \stackrel{i-th}{1}, \ldots, 0 ), \; 1 \leq i \leq 
\delta.$$On the other hand, since $\E$ is globally generated on $X$, the evaluation morphism
$$H^0(X, {\E}) \otimes {\Oc}_X \stackrel{ev}{\to} \E$$is surjective. This means 
that, for each $p \in X$, there exist global sections $s_1^{(p)}, \; s_2^{(p)} \in 
H^0 (X, {\E})$ such that 
$$s_1^{(p)}(p) = (1,0), \; s_2^{(p)}(p) = (0,1) \in {\Oc}_{X,p}^{\oplus 2}.$$Therefore, 
it immediately follows that 
$$H^0(X, {\E} \otimes M^{\otimes k}) \to \!\!\! \to H^0({\Oc}_{\Sigma}^{\oplus 2}) \cong 
{\CC}^{2 \delta}.$$If we compose with diagram \eqref{eq:(*1)}, we get: 
\begin{displaymath}
\begin{array}{ccl}
H^0( {\E} \otimes M^{\otimes k})  & \to\!\!\! \to & H^0({\Oc}_{\Sigma}^{\oplus 2}) \cong {\CC}^{2 \delta} \\ 
\downarrow^{\mu_X} & & \downarrow  \\ 
H^0({\Oc}_{\Sigma}) & \stackrel{\cong}{\to} & H^0({\Oc}_{\Sigma}) \cong  {\CC}^{\delta}.\\ 
 & & \downarrow \\ 
  & & 0 
\end{array}
\end{displaymath}thus $\mu_X$ is surjective (cf. Remark \ref{rem:localdescr}). 
By \eqref{eq:regbis}, one can conclude.
\end{proof}

\begin{remark}\label{rem:10}
\normalfont{
Observe that the bound (\ref{eq:uniform}) is uniform, i.e. it does not depend 
on the postulation of nodes of the curves which are zero-loci of sections parametrized
by ${\V}_{\delta}({\E} \otimes M^{\otimes k})$. We remark that Theorem \ref{thm:9bis} improves 
our Theorem 4.5 in \cite{F2}. Both these results generalize 
what proved by Ballico and Chiantini in \cite{BC} mainly 
because our approach more generally holds for families of nodal curves on smooth 
projective threefolds but also because, even in the case of $X = \Pt$, main subject 
of \cite{BC}, our regularity results are effective and not asymptotic as Proposition 3.1 in \cite{BC}. 
Furthermore, in \cite{F2} we observed that 
the bound $\delta \leq k+1$ is almost sharp. Indeed, one can easily 
construct examples of non-regular points 
$[s] \in  {\V}_{k+4}({\Oc}_{\Pt}(k+1)\oplus {\Oc}_{\Pt}(k+4))$, for any $ k \geq 3$, 
whose corresponding curve $C$ has its $(k+4)$ nodes lying on a line $L \subset \Pt$; anyhow, one 
can also show that 
${\V}_{k+4}({\Oc}_{\Pt}(k+1)\oplus {\Oc}_{\Pt}(k+4))$ is generically regular.  

} 
\end{remark}

For what concerns $\Pmc_{\delta}$-Severi varieties on $\Pmc$, we get:
\begin{theorem}\label{thm:regrdelta}
Let $X$ be a smooth projective threefold, $\E$ be a globally generated rank-two vector bundle on 
$X$, $M$ be a very ample line bundle on $X$ and $k \geq 0$ and $\delta >0$ be integers.

\noindent
Let $\Pmc : = \Pp_X (\E \otimes M^{\otimes k})$ and let $\Oc_{\Pmc}(1)$ be its tautological 
line bundle. Let ${\mathcal R}_{\delta} (\Oc_{\Pmc}(1))$ be the $\Pmc_{\delta}$-Severi 
variety of irreducible divisors on $\Pmc$ having $\delta$-rational double points on $\Pmc$. 
Then, if: 
\begin{equation}\label{eq:uniform2}
\delta \leq k+1,
\end{equation} ${\mathcal R}_{\delta} (\Oc_{\Pmc}(1))$  is regular.
\end{theorem}
\begin{proof} From Theorem \ref{thm:9bis}, we know that \eqref{eq:uniform2} is a sufficient 
condition for the regularity of  
${\V}_{\delta}({\E} \otimes M^{\otimes k})$ on $X$. One can conclude 
by using Corollary \ref{cor:tangentspace}. 
\end{proof}

\end{document}